\documentclass[10pt,a4paper,reqno]{amsart} 
\usepackage{latexsym,floatflt}
\usepackage{epsfig}
\usepackage{rotating}
\usepackage{amssymb}
\usepackage[T1]{fontenc}
\usepackage{afterpage}
\usepackage{color}

\catcode`\@=11
\def\section{\@startsection{section}{1}%
 \z@{.7\linespacing\@plus\linespacing}{.5\linespacing}%
 {\normalfont\bfseries\scshape\centering}}
\def\subsection{\@startsection{subsection}{2}%
  \z@{.5\linespacing\@plus\linespacing}{.5\linespacing}%
  {\normalfont\bfseries\scshape}}

\def\subsubsection{\@startsection{subsubsection}{3}%
 \z@{.5\linespacing\@plus\linespacing}{-.5em}%{.5\linespacing}%
  {\normalfont\bfseries\itshape}}
\catcode`\@=12

\def\cqfd{$\hfill{\vrule height 3pt width 5pt depth 2pt}$}

\newcommand{\ov}{\overline}

\newcommand{\f}{{\hbox{\small {\rm fr}}}}
\newcommand{\GK}{\mathbb{K}}
\newcommand{\GL}{\mathbb{L}}
\newcommand{\GM}{\mathbb{M}}
\newcommand{\ns}{\mathbb{N}}
\newcommand{\zs}{\mathbb{Z}}

\newcommand{\qs}{\mathbb{Q}}
\newcommand{\rs}{\mathbb{R}}

\newcommand{\cs}{\mathbb{C}}

\newcommand{\bm}[1]{\mbox{\boldmath \ensuremath{#1}}}

\newcommand{\R}{\mathcal R}

\newcommand{\Ref}[1]{(\ref{#1})}
\newcommand{\beq}{\begin{equation}}
\newcommand{\eeq}{\end{equation}}
\newcommand{\gf}{generating function}
\newcommand{\gfs}{generating functions}

\newcommand{\eps}{\epsilon}

\newcommand{\fps}{formal power series}
\def\cqfd{\par\nopagebreak\rightline{\vrule height 4pt width 5pt depth 1pt}
\medbreak}
\def\emm#1,{{\em #1}}

\newtheorem{Theorem}{Theorem}
\newtheorem{Propo}[Theorem]{Proposition}
\newtheorem{Coro}[Theorem]{Corollary}
\newtheorem{Lemma}[Theorem]{Lemma}
\newtheorem{Definition}[Theorem]{Definition}

\begin{document}

%%%%%%%%%%%%%%%%%%%%%%%%%%%%%%%%%%%%%%%%%%%%%%%%%%%%%%%%%%%%
\title[Polynomial equations with one catalytic variable (Section~\thesection)]
{ Polynomial equations with one catalytic variable,\\   algebraic series,\\
and  map enumeration }

\author{Mireille Bousquet-M\'elou \and Arnaud Jehanne}
\address{CNRS, LaBRI, Universit\'e Bordeaux 1, 351 cours de la Lib\'eration,
  33405 Talence Cedex, France \\
\and A2X, Universit\'e Bordeaux 1, 351 cours de la Lib\'eration,
  33405 Talence Cedex, France}
\email{mireille.bousquet@labri.fr, jehanne@math.u-bordeaux1.fr}
\thanks{MBM  was partially supported by the European Commission's IHRP
  Programme, grant HPRN-CT-2001-00272, ``Algebraic Combinatorics in
  Europe''} 
%\subjclass{Primary: subject; Secondary: subject}
\date{April 1, 2005}

%\date{}

\maketitle

\begin{abstract}
Let $F(t,u)\equiv F(u)$ be a \fps\ in $t$ with polynomial coefficients
in $u$. Let $F_1 , \ldots, F_k$ be $k$ \fps \ in $t$, independent of
$u$. Assume all these series are characterized by a polynomial
equation
$$
P(F(u), F_1, \ldots , F_k, t , u)=0.
$$
We prove that, under a mild hypothesis on the form of  this equation,
these $(k+1)$ series are 
algebraic, and we  give a strategy to compute a polynomial equation for
each of them. \\
This strategy generalizes the so-called \emm kernel method, and  \emm
quadratic method,, which apply respectively to equations that are
linear and quadratic in $F(u)$. Applications include the solution of
numerous map enumeration 
problems, among which the hard-particle model on general planar maps.
\end{abstract}

%%%%%%%%%%%%%%%%%%%%%%%%%%%%%%%%%%%%%%%%%%%%%%%%%
\section{Introduction}
\label{section-intro}
%%%%%%%%%%%%%%%%%%%%%%%%%%%%%%%%%%%%%%%%%%%%%%%%%%
Let us begin with a classical enumeration problem. We consider walks
on the half-line $\ns$, that start from $0$ and consist 
of unit steps $\pm 1$. Let $F(t,u) \equiv F(u)$ be their \gf , where
$t$ counts the 
length (the number of steps) and $u$ the position of the
endpoint. That is to say, $F(t,u) = \sum_{n,k} a_{n,k} t^n u^k$, where
$a_{n,k}$ is the number of $n$-step walks that end at level $k$. 
Note that $F(t,0) \equiv F(0)$ 
is the length \gf \ of the celebrated \emm Dyck
paths,, which are the walks ending at
$0$~\cite[p.~173]{stanley-vol2}. A step-by-step  
construction of these walks gives either a recurrence relation of the
numbers $a_{n,k}$ or, equivalently, the following functional equation:
\beq
\label{walk-1D}
F(u) = 1+ tuF(u) + \frac t u \Big( F(u)-F(0)\Big).
\eeq
The second (resp. third) term on the right-hand side counts walks
ending with a step $+1$ (resp. $-1$). Clearly, this equation defines
$F(u)$ uniquely as a \fps\ in $t$ (with 
rational coefficients in $u$). Observe that the equation
\beq
\label{walk-1D-F1}
F(u) = 1+ tuF(u) + \frac t u \Big( F(u)-F_1\Big)
\eeq
defines uniquely \emm both,\, $F(u)$ and $F_1$ as \fps\ in $t$, if we
impose that $F(u)$ has \emm polynomial,\, coefficients in $u$ and that
$F_1$ is independent of $u$. Indeed,
after multiplying the equation by $u$ and setting $u=0$, we find
$F_1=F(0)$, and we are thus back to~\Ref{walk-1D}. Finally, we recall
that $F(0)$ is well-known to be \emm algebraic,\, of degree $2$,
$$
F(0)=  \frac{1-\sqrt{1-4t^2}}{2t^2}.
$$
Consequently, $F(u)$ is algebraic too (meaning that it satisfies a
non-trivial polynomial equation, $Q(t,u,F(u))=0$, with rational
coefficients). 

\smallskip
The above example is an instance 
of the general situation we study in this paper. We assume that a
 $(k+1)$-tuple  $(F(u), F_1, F_2, \ldots ,
F_k)$ of  \fps \ in $t$ is completely determined by a polynomial equation
\beq
\label{eq-generic}
P(F(u), F_1, F_2, \ldots , F_k, t, u)=0.
\eeq
Typically, $F(u)$ has polynomial coefficients in $u$, and $F_i$
is the coefficient of $u^{i-1}$ in $F(u)$. Following Zeilberger's
terminology~\cite[p.~457]{zeil-umbral}, we say 
that~\eqref{eq-generic} is a \emm polynomial
equation with one catalytic variable, $u$. The aim of this paper
is twofold: we prove that the solution of a (well-founded) equation
of the form~\Ref{eq-generic} is \emm always algebraic,\, and we
present a strategy to obtain a polynomial equation it satisfies.

There are several reasons why we like to know that the \gf\ of some
class of objects is algebraic. Firstly, the set of algebraic series is closed
under natural operations (sum, product, derivatives,
composition...). Secondly, these series are reasonably easy to handle (via
resultants or Gr\"obner bases). In particular,  several computer algebra
packages are now able to make the above closure
properties effective. Thirdly, algebraic series are also \emm
D-finite, and this  
implies that their coefficients can be computed in a linear number of
operations~\cite[Ch.~6]{stanley-vol2}. The asymptotic behaviour of
these coefficients has 
a generic form, the details of which are usually not to hard to
obtain. Finally -- and, to many combinatorialists, most importantly --
the fact that a class of objects is counted by an algebraic series
suggests that it should be possible to construct these objects
recursively by \emm concatenation, of objects of the same type. For
many objects, such a construction is easily found, but for others,
among which planar maps~\cite{cori-vauquelin,Sch97}, the algebraic structure of the objects is far
from clear, and the algebraicity of the \gf\ gives rise to challenging
combinatorial problems. See~\cite[Ch.~6]{stanley-vol2}
or~\cite{fla-poly-alg} for a presentation of algebraic series in
enumeration. 

%\medskip

But let us return to polynomial equations with one catalytic variable.
Many combinatorialists have fought them before us, and we
want to recall some milestones in this history.

\subsection{A partial historical account}
In 1956 already, Temperley writes, for the perimeter enumeration of
\emm column-convex polyominoes,\,, a set of recurrence
relations~\cite[Eq.~(7)]{temperley} that 
is equivalent, after summation, to
$$
F(u)= \frac{ut^2}{1-ut} +\frac{t^3u^2 F(u)}{(1-ut)^2}
+2\,\frac{t^2u^2}{1-ut}\,\frac{F(u)-F(1)}{u-1}
+ut\, \frac {uF(u)-uF(1)-(u-1)F'_u(1)}{(u-1)^2}.
$$
He proves that $F(1)$ is algebraic, without being able to compute it
explicitly (see~\cite{feretic} for a simple expression of
$F(1)$). Like~\eqref{walk-1D}, the above equation is linear in $F(u)$, but it
contains \emm two, additional unknown functions, $F(1)$ and $F'_u(1)$.

The first non-linear equations appear in the early sixties, in the
work of Tutte and Brown on planar maps. For instance, Tutte publishes
in 1962 the following equation~\cite[Eq.~(3.7)]{tutte-triangulations},
which rules the enumeration of certain triangulations:
\beq
\label{eq-triangulations}
F(u)=1 + \frac t u \, \left( \frac{F(u)}{1-uF(u)} -F(0)\right).
\eeq
In the following years, more equations of this type are
published for various families of planar maps
(non-separable~\cite{brown-non-separable,brown-tutte-non-separable}, 
general~\cite{tutte-general}, other triangulations~\cite{brown-triangulations},
quadrangulations~\cite{brown-quadrangulations}). All of them 
involve only one unknown 
function $F_1$ (and thus read $P(F(u),F_1,t,u)=0$),  and are quadratic
in $F(u)$, apart from the equation on 
quadrangulations which is cubic. In the first papers, Tutte and Brown
solve these 
equations by \emm guessing and checking,\,: either they  guess the
expansions of $F_1$ and $F(u)$, and
then check that their guesses satisfy the functional equation, or they
only guess the expansion of $F_1$, and then prove that the polynomial
equation $P(F(u),F_1,t,u)$, taken with the conjectured value of $F_1$,
admits one root $F(u)$ that is a \fps\ in $t$ with polynomial
coefficients in $u$. 
Of course, any equation for which the value of $F_1$ cannot be guessed
remains hopeless with this strategy.  

\medskip
In 1965, Brown publishes a theorem that deals, at first
sight, with a different topic: with the conditions satisfied by a series
in $t$ and $u$ that admits a square root (which is itself assumed to be
a \fps)~\cite{brown-square}. He
shows that this theorem allows to solve, in a systematic way, all
equations of the type~\Ref{eq-generic} that are quadratic in $F(u)$
and only involve one unknown function $F_1$. The
\emm quadratic method,\,  is born
(see~\cite[Section~2.9]{goulden-jackson} for a modern 
account). Brown even manages to solve, with some contorsions, 
the above-mentioned cubic equation for
quadrangulations~\cite[Section~4]{brown-square}. At the end 
of~\cite{brown-conf}, he writes \emm ``It is possible that the method
may be effective when more than one unknown series is present'',. This
hope was  confirmed many years later, in
1994, when Bender and Canfield applied the method to a quadratic
equation with arbitrarily many  unknown 
functions~\cite{bender-canfield}.

But let us go back to the sixties. In 1968, in the first volume
of \emm The art of computer programming,\,, Knuth gives for the
classical \emm ballot problem, an equation that is equivalent
to~\Ref{walk-1D}, and presents a ``trick'' that solves
it~\cite[Section~2.2.1, Ex.~4]{knuth}. This may have been 
the unnoticed birth of the \emm kernel 
method,\,, which allows to solve systematically equations of the
form~\Ref{eq-generic} that are linear in $F(u)$. This trick may
have been better known at that time in probability theory. At least,
the same idea definitely appears in a 1979
paper~\cite{fayolle-processors}, in a more difficult, analytic
context. The kernel method is currently the subject of a certain
revival in
combinatorics~\cite{hexacephale,banderier-flajolet,bousquet-petkovsek-1,de-mier,prodinger}.

In 1972, Cori and Richard solve again certain linear equations, and
also some polynomial equations  with one unknown series
$F_1$~\cite{cori-richard}. Their technique 
is very interesting, but the fact that they
deal with equations in  non-commuting variables makes it both
deeper and more obscure. Still, the strategy we present here to
attack~\Ref{eq-generic} owes a lot to~\cite{cori-richard}.

\medskip
Since then, equations of the form~\Ref{eq-generic} have continued to
appear in various enumeration problems, mostly involving
maps~\cite{bender-canfield,gao-5-connected}, but also
polyominoes~\cite{bousquet-habilitation,bousquet-vcd,feretic-svrtan},
stack-sortable 
permutations~\cite{bousquet-stack,zeilberger-stack} and their
generalizations~\cite{xu}, lattice
walks~\cite{banderier-flajolet,bousquet-versailles}, etc. 
Examples can be found were both the degree of the equation and the
number of unknown functions is arbitrarily large. For instance, such
equations are hiding in Tutte's work on the chromatic polynomial of
triangulations~\cite[Section~5]{tutte-chromaticV}.
 Another such set of equations, unbounded in degree and
number of unknowns, is presented  in
Section~\ref{subsection-constellations}. It deals with the enumeration of
certain Eulerian maps called \emm constellations,.

\subsection{Contents}
\noindent {\bf \em The general strategy.}
The method presented in this paper to solve equations of the
form~\Ref{eq-generic} encapsulates and simplifies all previous
approaches, in particular the kernel method and the quadratic
method. It works without any restriction on the degree of the equation
or on the number of unknowns $F_i$. The general strategy is described
in Section~\ref{section-key}. Its justification only takes a few
lines. It yields a system of $3\ell$ polynomial equations that relate
$k+2\ell$  series in $t$: the unknowns $F_1,\ldots , F_k$, first,
then $\ell$ series named $U_1, \ldots , U_\ell$, which are defined as the
roots of a certain equation (simply related to the original functional
equation), and finally the values 
of $F(u)$ at $u=U_i$, for $i=1, \ldots , \ell$. The strategy ``works''
if, first, 
$\ell=k$ (so that we have as many equations as unknown series), and if
 the $3k$ polynomial equations thus obtained imply the
algebraicity of the $F_i$. 

\smallskip
\noindent {\bf \em First examples.}
In Section~\ref{section-examples}, we apply this strategy to several
examples. For each of them, we observe that the strategy 
works: we find as many series $U_i$ as we have unknowns
$F_i$, and we can derive from the system of $3k$ polynomial equations
 an algebraic equation for each $F_i$.   We also
 relate our approach to the earlier \emm kernel method, and \emm
 quadratic method,.

\smallskip
\noindent {\bf \em A generic algebraicity theorem.}  Will this 
strategy \emm always, work? 
Section~\ref{section-generic} answers
this question positively, at least for  a well-founded equation of
the form
\beq
\label{generic-Q}
F(u) = F_0(u)+ t\  Q\Big( F(u), \Delta F(u), \Delta ^{(2)}F(u),\ldots
,\Delta ^{(k)}F(u), t, u\Big) ,
\eeq
where $F_0(u)$ is a given polynomial in $u$ (with 
coefficients in a field $\GK$ of characteristic $0$), $Q(x_0, \ldots ,
x_k,t,v)$ is another polynomial, and 
$$
\Delta ^{(i)} F(u) = \frac{F(u)-F_1-uF_2 -\cdots -
  u^{i-1}F_i}{u^i}.
$$
If we require $F(u)$ to be a \fps\ in $t$ with polynomial coefficients
in $u$, and  $F_i$ to be  the coefficient of
$u^{i-1}$ in $F(u)$, then  this equation defines uniquely
$F(u)$.  We prove that $F(u)$,  and hence all the $F_i$, are
algebraic. 

\smallskip
\noindent{\bf \em Algebraicity results for planar maps.}  Thus the solution of
every (well-founded) equation with one catalytic variable  is
algebraic. This result urges a combinatorial interlude, in which
we establish for several families of planar maps an equation of this
type. The generic algebraicity theorem tells us, without going
further, that their \gfs\ are  algebraic. Our
examples include some already studied problems (like the
face-distribution of Eulerian maps, for which we answer positively a
question left open in~\cite{mbm-schaeffer-ising}), and some new ones, like the
hard-particle model on general planar maps.

\smallskip
\noindent{\bf \em From $3k$ to $2k$, and then $k$ equations.}
The next question that we address is both theoretical and practical:
it deals with the size of our 
polynomial system. Assume the general strategy works and
provides a system of $3k$ equations. Even 
when $k=2$, even for a computer algebra system, this can be hard to
handle. In Section~\ref{section-discriminant} we reduce the system to
$2k$ equations which involve only the series $F_i$ and $U_i$. This new
system can be described simply in terms of the \emm discriminant,\, of
the polynomial $P$ occurring in~\Ref{eq-generic},
taken with respect to 
its first variable (we assume that $P$ is at least quadratic in
this variable). 
 Our $2k$ equations say that  
this discriminant, evaluated at $F_1, \ldots , F_k, t, u$ and considered
as a polynomial in $u$, 
\emm has $k$ multiple roots, $U_1, \ldots , U_k$.
This extends a result that was known to hold in the
quadratic case and is one of the possible formulations of the
quadratic method~\cite[Section~2.9]{goulden-jackson}.

%Another formulation of the above statement is that 
Hence the discriminant and its derivative with respect to $u$
have $k$ roots in common.
It is well-known that two polynomials have \emm one,\, root in common if their
resultant is zero. In Section~\ref{section-resultants}, we recall how to
express, by a set of $k$ determinants, the fact that two polynomials have $k$
roots in common. Applying this to the discriminant
and its derivative, we obtain a set of $k$ polynomial equations that relate
$F_1, \ldots , F_k$.

\smallskip
\noindent{\bf \em A new proof of Brown's theorem.}
Before turning our attention to specific examples, we  give  in
  Section~\ref{section-brown} a
``modern'', and maybe clearer proof of Brown's theorem on square
roots of bivariate power series\footnote{As mentioned in
  Section~\ref{section-brown}, it seems that there may be
  a mistake in Brown's original proof.}.
Recall that this theorem is the basis of the quadratic method. 

\smallskip
\noindent{\bf \em Practical examples.} We discuss in
Section~\ref{section-practise} 
how to derive in practise  an
algebraic equation for, say, the unknown series $F_1$. We suggest
various approaches, which we exemplify on certain maps called
$3$-constellations. The associated equation is cubic and involves two
unknown series $F_i$. In Section~\ref{section-distribution}, we walk
in the steps of 
Bender and Canfield~\cite{bender-canfield} to find the
face-distribution of planar maps. This 
problem was already solved in two other ways~\cite{BDG-planaires}, and
we prove that our results are equivalent to the former
ones. Finally,
we solve in Section~\ref{section-hard} the hard-particle model on
general planar maps. For other recent applications of our method,
see~\cite{bernardi}.

\smallskip
Finally, Section~\ref{section-final} discusses a number of open questions.

%%%%%%%%%%%%%%%%%%%%%%%%%%%%%%%%%%%%%%%%%%%
\subsection{Formal power series and their relatives}
Let us conclude  this introduction with some  notation.
%  definitions and of formal
%series. \\
Let $\GK$ be a commutative ring. We denote by
$\GK[t]$ the set of polynomials in $t$ with coefficients in
$\GK$. If $\GK$ is a field, then $\GK(t)$ denotes the field of
fractions in $t$ with coefficients in $\GK$. We denote by $\overline
\GK$ the algebraic closure of $\GK$.
We also consider several sets of series of the form
$$
A(t)=\sum_{n\ge n_0} a_n t^{n/d},
$$
where $n_0 \in  \zs$, $a_{n_0}\not = 0$  and $d\in
\ns\setminus\{0\}$. The number $n_0/d$ is called the \emm valuation,\,
of $A(t)$. We use the standard notation for the coefficients of a
series:
$$
[t^{n/d}]A(t):= a_n.
$$
In particular, 
\begin{itemize}
%-- $\GK(t)$ the set of rational functions of  with coefficients in
%$\GK$,\\
\item[--]  $\GK[[t]]$ is the set of formal power series in $t$ with coefficients
in $\GK$ ($n_0\ge 0$ and $d=1$),
\item[--]  $\GK((t))$ is the set of Laurent series in $t$ with coefficients
in $\GK$ ($d=1$),
\item[--]  $\GK^\f[[t]]$ is the set of fractional power series in $t$
 with coefficients in $\GK$ ($n_0\ge 0$),
\item[--]  $\GK^\f((t))$ is the set of fractional Laurent  series in
  $t$ (a.k.a. Puiseux series) with coefficients in $\GK$ (no condition).
\end{itemize}
Each of these sets is a commutative ring, and the second and fourth
are  fields
if $\GK$  is a field. More precisely, $\GK((t))$ is the fraction field
of $\GK[[t]]$, and  $\GK^\f((t))$ is the fraction field
of $\GK^\f[[t]]$. If, moreover, $\GK$ is 
 algebraically closed and has characteristic $0$, then so is
 $\GK^\f((t))$~\cite[Thm.~6.1.5]{stanley-vol2}. 

These notations  generalize to series in several
indeterminates. In this paper, we will mostly use series in $t$ and
$u$. Note the following inclusions:
$$
\GK[t,u] \subset \GK[[t]][u] \subset \GK[u][[t]] \subset \GK[[t,u]]=
\GK[[t]][[u]]. 
$$
The second set above is the set of polynomials in $u$ whose
coefficients are \fps\ in $t$. The third set is the set of \fps\ in $t$
whose coefficients are polynomials in $u$. The notation
$\GK[[u]]^\f[[t]]$ stands for the set of power series in $u$ and $t$
that are fractional in $t$.

All the fields considered in this paper have
implicitly 
characteristic $0$.,

%%%%%%%%%%%%%%%%%%%%%%%%%%%%%%%%%%%%%%%%%%%%%%%%%
\section{The general strategy}\label{section-key}
%%%%%%%%%%%%%%%%%%%%%%%%%%%%%%%%%%%%%%%%%%%%%%%%%%
Let $\GK$ be a field.
In our examples, $\GK$
will be $\cs$, or a field of fractions like $\cs(s_1,\ldots, s_m)$.
Let $F(t,u)\equiv F(u)$ be a series of $\GK[u][[t]]$, and let $F_1(t)
\equiv F_1, \ldots , F_k(t)\equiv F_k$ be $k$ series of $\GK[[t]]$. In
our framework, these $k+1$ series are the \gfs\ of certain families of
objects, counted according to one or two parameters. Assume these
series are related by an equation of the form
\beq
\label{eq-generic1}
P(F(u), F_1, F_2, \ldots , F_k, t, u)=0,
\eeq
where $P(x_0,x_1, \ldots , x_k, t, v)$ is a non-trivial polynomial in
$k+3$ variables, with coefficients in $\GK$. Assume, moreover, that
the above equation defines the $(k+1)$-tuple $(F(u),F_1, \ldots , F_k)$
\emm uniquely,\, in the set $\GK[u][[t]] \times \GK[[t]]^k$. 
Some examples were given in the introduction, and numerous examples
will be given below. 

Let us differentiate~\Ref{eq-generic1} with respect to $u$:
$$
F'(u) 
\frac {\partial P}{\partial x_0} (F(u),F_1, \ldots , F_k, t,u)
+\frac {\partial P}{\partial v} (F(u),F_1, \ldots , F_k, t,u)=0.
$$
Let $U(t)\equiv U$ be a series of $\GK^\f[[t]]$. The series
$F(U)\equiv F(t,U)$ is a well-defined fractional power series in
$t$. The same holds for $F'(U)$. If,
moreover,
\beq
\label{eq-derivee-x0}
\frac {\partial P}{\partial x_0} (F(U),F_1, \ldots , F_k, t,U)=0,
\eeq
then the above identity implies that 
$$
\frac {\partial P}{\partial v} (F(U),F_1, \ldots , F_k, t,U)=0.
$$

This simple observation is the key of our solution of equations of the
form~\Ref{eq-generic1}. {\em If we can prove the existence of $k$
distinct series $U_1, \ldots , U_k$, belonging to  $\GK^\f[[t]]$, that
satisfy~{\em\Ref{eq-derivee-x0}}, then the following system of $3k$
polynomial equations holds}: for $1\le i\le k$,
\begin{eqnarray}
P\Big(F(U_i),F_1, \ldots , F_k, t,U_i\Big) &=&0,
\label{key-P}\\
\frac {\partial P}{\partial x_0} \Big(F(U_i),F_1, \ldots , F_k, t,U_i\Big)
&=&0,
\label{key-PF}\\
\frac {\partial P}{\partial v}\Big(F(U_i),F_1, \ldots , F_k, t,U_i\Big)
&=&0.\label{key-Pu}
\end{eqnarray}
A bit of optimism allows us to hope that this system characterizes
completely the $3k$ unknown series it involves, namely $F_1, \ldots ,
F_k$, $U_1, \ldots , U_k$ and $F(U_1), \ldots , F(U_k)$, so that each
unknown series (in particular each $F_i$) is algebraic. More precisely, we
would like this 
system to have only a finite number of solutions \emm under the assumption
that the series $U_i$ are distinct,. 
%If this is the case, it remains
%to compute the first few coefficients of the $F_i$ 
 This assumption can be encoded by adding a new unknown
$X$ and a new polynomial equation:
\beq
X\prod_{1\le i< j\le k }(U_i-U_j)=1.
\label{Ui-distinct}
\eeq

We prove in Section~\ref{section-generic} that this optimism is
justified: the solution of a well-founded equation of the
form~\Ref{generic-Q} is indeed shown to be algebraic. However, we do
not need this general theorem to examine and solve \emm specific, examples,
like~\Ref{walk-1D-F1} or~\Ref{eq-triangulations}. What we \emm do,\,
need is a way to determine 
how many series $U$ satisfy~\Ref{eq-derivee-x0}, without knowing the
value of $F(u)$ or $F_1, \ldots , F_k$. This turns out to be easy.
 Let us first clarify what we mean by a \emm root, \ $U$ of a
series $\Phi(t,u)$.
\begin{Lemma}
\label{lemma-factorization}
Let $\Phi(t,u)\in \GK[u]^\f[[t]]$, and $U\in \GK^\f[[t]]$. Then
$\Phi(t,U)$ is a well-defined series of $\GK^\f[[t]]$. If this series
is zero, we say that $U$ is a root of $\Phi(t,u)$. In this case, there
exists $\Psi(t,u)\in \GK[u]^\f[[t]]$ such that
$$
\Phi(t,u)=(u-U)\Psi(t,u).
$$
More generally, if $\Phi(t,u)$ factors as
$$
%\Phi(t,u)=(u-U_1)^{m_1}\cdots (u-U_k)^{m_k}\Psi(t,u),
\Phi(t,u)=(u-U)^{m}\Psi(t,u),
$$
where $\Psi(t,u)\in \GK[u]^\f[[t]]$, the series $U$ belongs to
$\GK^\f[[t]]$
% and are distinct, 
and $\Psi(t,U)\not =0$, we say that $U$ is a root of  $\Phi(t,u)$ of
multiplicity $m$.

This extends to the case where  $\Phi(t,u)$ belongs to
$\GK[[u]]^\f[[t]]$, if we require that $U$ has no constant term (that
is, vanishes at $t=0$). In this case,  $\Psi(t,u)$ also belongs to
$\GK[[u]]^\f[[t]]$.
\end{Lemma}
\noindent
{\bf Proof.}
The fact that $\Phi(t,U)$ is  well-defined is obvious, by definition
of the substitution of series: If
$$
\Phi(t,u)=\sum_{n\ge 0} t^{n/d} \phi_n(u),
$$
where $\phi_n(u)$ is a polynomial in $u$, then
$$
\Phi(t,U)=\sum_{n\ge 0} t^{n/d} \phi_n(U),
$$
and the coefficient of $t^{p/q}$  in $\Phi(t,U)$, for $p/q\le{k/d}$,
depends only on the polynomials $\phi_0(u), \ldots , \phi_k(u)$.
Now for any indeterminate $v$,
$$
\Phi(t,u)-\Phi(t,v)=(u-v)\sum_{n\ge 0} t^{n/d} \phi '_n(u,v),
$$
where 
$$
\phi '_n(u,v)= \frac{\phi_n(u)-\phi_n(v)}{u-v}
$$
is a polynomial in $u$ and $v$. The case $v=U$  proves the second
statement of the lemma. 

The argument can be adapted without any difficulty to the case where
$\Phi(t,u)$ belongs to 
$\GK[[u]]^\f[[t]]$ and $U$ has no constant term, upon writing
$$
\Phi(t,u) =\sum_{n, m \ge 0} \phi_{m,n} u^m t^{n/d} 
$$
with $\phi_{m,n}\in \GK$.
\cqfd

The next theorem tells how many
roots a series $\Phi(t,u)$ has.
\begin{Theorem}\label{thm-roots}
Let $\Phi(t,u) \in \GK[u]^\f[[t]]$, where $\GK$ is an algebraically
closed field.
% of characteristic $0$. 
Assume that the coefficient of $t^0$
in $\Phi$, that is to say, the polynomial $\Phi(0,u)$, is non-zero and
has degree $k$. Then $\Phi(t,u)$ has exactly $k$ roots in
$\GK^\f[[t]]$, counted with multiplicities. Let $U_1, \ldots , U_k$
denote these roots. Then
$$
\Phi(t,u) = (u-U_1) \cdots (u-U_k) \Psi(t,u)
$$
where $\Psi(t,u) \in \GK[u]^\f[[t]]$.
\end{Theorem}
\noindent
{\bf Proof.} The proof is a harmless extension of the proof of the Puiseux
theorem, which establishes the above result (and more) in the case where
$\Phi(t,u) \in \GK^\f((t))[u]$. We refer the reader
to~\cite[Ch.~4]{walker}. 
The coefficients of the $U_i$ can be computed
inductively using Newton's polygon.
\cqfd

\noindent

%%%%%%%%%%%%%%%%%%%%%%%%%%%%%%%%%%%%%%%%%%%%%%%%%
\section{First examples}
\label{section-examples}
%%%%%%%%%%%%%%%%%%%%%%%%%%%%%%%%%%%%%%%%%%%%%%%%%
We now apply our general strategy to a few examples.
\subsection{Walks on a half-line and the kernel method}
\label{section-kernel}
%\end{Example}
We consider here some equations of the
form~\Ref{eq-generic1} that are linear in $F(u)$. The reader familiar
with the \emm kernel method, will not find our calculations very original,
and this is normal: beyond solving these equations,  our
objective here is to show that our general strategy reduces to the
kernel method when the equation is linear. We refer
to~\cite{banderier-flajolet} for a systematic treatment of walks on
the half-line, based on the kernel method.

Let us first go back to the simplest equation we have met so far,
Eq.~\Ref{walk-1D-F1}. 
%We consider walks on the half-line $\ns$, that
% start at $0$ and consist 
%of steps $\pm 1$. Let $F(t,u) \equiv F(u)$ be their \gf , where $t$ counts the
%length (the number of steps) and $u$ the endpoint. Note that $F(t,0)
% \equiv F(0)$ 
%is the length \gf \ of the celebrated \emm Dyck walks,. A step-by-step
%construction of these walks gives
%$$
%F(t,u) = 1+ tuF(u) + \frac t u \Big( F(u)-F(0)\Big).
%$$
It  can be rewritten under the form~\Ref{eq-generic1}:
$$
%\Big(u-t(1+u^2)\Big)F(u)-u+tF(0)=0.
P(F(u),F_1, t,u)=0,
$$
where
$$
P(x_0,x_1,t,v)=\left(v-t\left(1+v^2\right)\right)x_0-v+tx_1.
$$
  Condition~\Ref{eq-derivee-x0} reads in this case:
$$
U-t(1+U^2) =0.
$$
In accordance with Theorem~\ref{thm-roots}, we find that there exists
a unique fractional power series in $t$ that satisfies this equation, namely
$$
U= \frac{1-\sqrt{1-4t^2}}{2t}.
$$
The system~(\ref{key-P}--\ref{key-Pu}) now reads
\begin{eqnarray*}
\left(U-t\left(1+U^2\right)\right)F(U)&=&U-tF_1,\\
U-t(1+U^2) &=&0,\\
(1-2tU) F(U) &=&1.
\end{eqnarray*}
The first and second equations together imply that
$$
F_1= \frac U t= \frac{1-\sqrt{1-4t^2}}{2t^2}.
$$
We have  recovered the classical expression of the \gf \ of Dyck
paths.
%If necessary, 
%the third equation gives $F(U)$, and 
An expression for $F(u)$ now follows from the original equation
$P(F(u),F_1, t,u)=0 $. 

\medskip
Let us now study a problem with more unknown functions. We still
consider walks on the half-line $\ns$ that start from $0$, but
they now consist of steps $+3$ and $-2$. A step-by-step construction of
these walks gives, for their bivariate \gf\ $F(t,u)\equiv F(u)$, the
equation
\beq
\label{+3-2}
F(u)=1+ tu^3F(u) +\frac t{u^2}\left( F(u)-F_1-uF_2\right)
\eeq
where $F_1$ (resp.~$F_2$) is the length \gf\ of walks ending at $0$
(resp.~$1$). This equation can be rewritten as
$
P(F(u), F_1,F_2,t,u)=0,
$
with
$$
P(x_0,x_1,x_2,t,v)= \left( {v}^{2}-t ( 1+{v}^{5})  \right)x_0
-{v}^{2} +t{x_1}+tv{ x_2}.
$$
 Condition ~\Ref{eq-derivee-x0} now reads
$$
U^2-t(1+U^5)=0.
$$
By Theorem~\ref{thm-roots}, exactly two fractional power
series $U_1$ and $U_2$ satisfy  this equation, and we happily observe
that two is also the number of unknown series $F_i$. 
%  The solutions of the above equation are actually formal power series 
% in $\sqrt t$, and 
One may compute the first terms of the $U_i$'s
using Newton's polygon:
$$
U_{1,2}=\pm  t^{1/2}+\frac 1 2\, t^3 \pm \frac 9 8\,  t^{11/2}+\frac 7
2\, t^8+O(t^{21/2}). 
$$
In particular, these two series are distinct. The
system~(\ref{key-P}--\ref{key-Pu}) now reads, for $i=1,2$,
\begin{eqnarray}
\left( {U_i}^{2}-t ( 1+{U_i}^{5})  \right)F(U_i)&=& {U_i}^{2} -t{F_1}-tU_i{
  F_2},
\label{complet2-3}\\
U_i^2-t(1+U_i^5)&=&0,\label{kernel2-3}\\
  U_i\left( 2-5\,t{U_i}^{3} \right)F(U_i) &=&2\,U_i-tF_2.\nonumber
\end{eqnarray}
We have thus obtained six equations that relate
$F_1,F_2,U_1,U_2,F(U_1)$ and $F(U_2)$. At this point, there are
several ways to conclude. The fastest one is probably to observe that,
by~\Ref{complet2-3} and~\Ref{kernel2-3}, the series $U_1$ and $U_2$
are the two roots of the following polynomial in $u$:
$$
R(u)=u^2-tuF_2-tF_1.
$$
Thus this polynomial factors as $(u-U_1)(u-U_2)$, which implies
$$
-tF_1=U_1U_2\quad \hbox{and} \quad tF_2=U_1+U_2.
$$
One can then eliminate $U_1$ and $U_2$ using~\Ref{kernel2-3}, and
obtain polynomial equations for $F_1$ and $F_2$. In particular, the
\gf\ $F_1$ of walks ending at $0$ satisfies:
$$
{F_1}=1+2\,{t}^{5}{{F_1}}^{5}-{t}^{5}{{F_1}}^{6}+{t}^{5}{{
F_1}}^{7}+{t}^{10}{{F_1}}^{10}.
$$

Consider, more generally, the case where the functional
equation~\Ref{eq-generic1} has degree~$1$ in $F(u)$ and can be written
as
$$
K(t,u) F(u) = P(F_1, \ldots , F_k,t,u)
$$
where $K(t,u) \in \GK[t,u]$ is the \emm kernel,\, of the equation, and 
$P(x_1, \ldots , x_k, t, u)$ is a polynomial in $k+2$ indeterminates.
The system~(\ref{key-P}--\ref{key-Pu}) reads
\begin{eqnarray*}
K(t,U) F(U)& =& P(F_1, \ldots , F_k,t,U),
\\
K(t,U) &=&0,
\\
K'_u(t,U) F(U)& =& P'_u(F_1, \ldots , F_k,t,U).
\end{eqnarray*}
By combining the first and second equations, we see that every root of
the kernel that is finite at $t=0$ gives a polynomial equation
relating the $k$ unknown series $F_1, \ldots , F_k$. This is exactly
the principle of the \emm kernel method,, which  has been
around since the  late
60's, and is currently the subject of a certain revival
(see~\cite{hexacephale,banderier-flajolet,bousquet-petkovsek-1,de-mier,prodinger}
and references therein).

\subsection{Planar maps and the quadratic method}
\label{section-quadratic}
We consider here \emm rooted planar maps, (see Section~\ref{section-equations}
 or~\cite{goulden-jackson} for  definitions). Let $F(t,u)\equiv
F(u)$ be their 
\gf , where $t$ counts the number of edges, and $u$ the degree of the
root-face. Deleting the root-edge gives~\cite[Eq.~(4)]{tutte-general}:
\beq
\label{eq-planar-maps}
F(u)  =1+tu^2F(u)^2+ tu \, \frac{uF(u)-F(1)}{u-1}.
\eeq
Multiplying this equation by $(u-1)$ gives a polynomial equation of
the form~\Ref{eq-generic1}, with one unknown function $F_1:=F(1)$.
 Condition~\Ref{eq-derivee-x0} reads in this case:
$$
U-1=2tU^2(U-1)F(U)+tU^2.
$$
By Theorem~\ref{thm-roots}, this equation has a (unique) solution $U$
in the set of fractional power series in $t$. It is actually clear
on the equation  that such a series exists, and is a  \fps\ in
$t$ (think of extracting the coefficient of $t^n$). Moreover, $U\not
=0,1$. From~(\ref{key-P}--\ref{key-Pu}), we obtain 
\begin{eqnarray*}
(U-1)F(U)&=&U-1+tU^2(U-1)F(U)^2+tU^2F(U)-tUF_1,\\
U-1&=&2tU^2(U-1)F(U)+tU^2,\\
F(U) &=& 1+ tU(3U-2)F(U)^2+2tUF(U)-tF_1.
\end{eqnarray*}
One can  eliminate $F(U)$ between the first and second equation, and
then between the second and the third. This gives two equations
relating $U$ and $F_1$. We ignore the irrelevant factors $U$ and
$U-1$, and eliminate $U$. This gives
an algebraic equation satisfied by $F_1$,  containing three distinct
factors. The right one is easily identified, given that $F_1=1+O(t)$,
and one concludes that the \gf\ of planar maps, counted by  edges, satisfies
$$
F_1=1-16\,t+18t F_1 -27\,{t}^{2}{F_1 }^{2}.
$$

\medskip
 More generally, an
equation of the form~\Ref{eq-generic1} having degree $2$ in $F(u)$ can
be written as
$$
\Big(2 a F(u)+b\Big)^2  = b^2-4ac = \Delta(u),
$$
where $a,b,c$ and $ \Delta$ lie in $\GK[t,u, F_1, \ldots , F_k]$.
The system~(\ref{key-P}--\ref{key-Pu}) reads
\begin{eqnarray*}
\Big(2 a F(U)+b\Big)^2   &=& \Delta(U),
\\
2 a F(U)+b &=&0,
\\
2\Big(2 a'_u F(U)+b'_u\Big)\Big(2 a F(U)+b\Big)   &=& \Delta'_u(U).
\end{eqnarray*}
By combining the first and second equations, we see that every
fractional power series $U$  that cancels
$2 a F(u)+b$  cancels the discriminant $\Delta$.
% and thus gives a polynomial equation relating the $k$ unknown series
% $F_1, \ldots , F_k$. 
By  combining the
second and third equations, we see that $U$ is actually a \emm multiple
root of the discriminant,. 

When there is only one unknown function
$F_1$, we recover exactly the \emm quadratic method,, as described
in~\cite{goulden-jackson}: if there exists a series $U$ such that $2 a
F(U)+b=0$, then $\Delta(u)$ admits  a multiple root. Hence the
discriminant of $\Delta(u)$ with respect to $u$ is zero: this gives an
algebraic equation satisfied by $F_1$.

This will be generalized in this paper to functional equations of the
form~\Ref{eq-generic1} and of degree at least two in $F(u)$: we will prove
that the 
discriminant  $\Delta$ of $P$, taken with respect to its first variable and
evaluated at $F_1, \ldots , F_k, t,u$, admits each $U_i$ as a multiple
root
% We will prove that any series $U$ that satisfies 
%$$
%\frac {\partial P}{\partial x_0} \Big(F(U),F_1, \ldots , F_k, t,U\Big)
%=0
%$$
%is a multiple root of $\Delta (P,x_0) (F_1, \ldots , F_k,t,u)$
(Section~\ref{section-discriminant}). 
%We also show in Section~\ref{section-resultants} how to express the
%fact that two polynomials have $k$ roots in common. 

\subsection{Quadrangular dissections of the disk}
\label{section-dissections}
Let us now consider a \emm cubic, example with one unknown function. This
example was solved by Brown, with some
difficulties~\cite{brown-square}.  Our 
strategy works without any restriction on the degree of the equation,
and the solution of this cubic example will be just as easy as the
solution of, say, the quadratic equation~\Ref{eq-planar-maps}.

The \emm quadrangular dissections of the disk, studied by Brown
in~\cite{brown-quadrangulations} can be described as the rooted,
non-separable planar maps, with no multiple edges, in which each
non-root face has degree 4 (see Section~\ref{section-equations} for
definitions). It is easy to see that the root-face of 
such maps has an even degree, at least equal to 4. 
% Also, the number of vertices is always even. 
Let $a_{n,k}$ be
the number of such maps with $n+4$ vertices in which the root-face has
degree $2k$,  and let
$$
F(t,u)\equiv F(u)= \sum_{n\ge0, k\ge 2} a_{n,k} t^n u^{k-2}.
$$
Eq.~(5.1) of~\cite{brown-quadrangulations} can be rewritten
as
$$
F(u)= \frac{F(u)-F_1}{u} -t^2F_1F(u) +2tF(u)(1+ut^2F(u))+(1+ut^2F(u))^3,
$$
where $F_1\equiv F(0)$ is 
% the coefficient of $u^0$ in $F(u)$, that is, 
the \gf\ of dissections of squares.
 Condition~\Ref{eq-derivee-x0} reads:
$$
U=
1-U{t}^{2}{ F_1}+2\,Ut \left( 1+2U{t}^{2}F(U) \right)
+3\,{U}^{2}{t}^{2} \left( 1+U{t}^{2}F(U) \right) ^{2}.
$$
By Theorem~\ref{thm-roots}, this equation has a (unique) solution $U$
in the set of fractional power series in $t$. (It is again clear
on the equation itself that such a series exists, and is a  \fps\ in
$t$.) Moreover, $U\not = 0$.  From~(\ref{key-P}--\ref{key-Pu}), we obtain 
\begin{eqnarray*}
UF(U)&= &{F(U)-F_1} -Ut^2F_1F(U) +2UtF(U)(1+Ut^2F(U))+U(1+Ut^2F(U))^3,\\
U&=&
1-U{t}^{2}{ F_1}+2\,Ut \left( 1+2U{t}^{2}F(U) \right) 
+3\,{U}^{2}{t}^{2} \left( 1+U{t}^{2}F(U) \right) ^{2},\\
F(U)&=&-{t}^{2}{F_1}\,F(U)+2\,tF(U) \left( 1+2U{t}^{2}F(U) \right)+ 
\left( 1+U{t}^{2}F(U) \right) ^{ 2}(1+4U{t}^{2}F(U)).
\end{eqnarray*}
One can  eliminate $F(U)$ between the first and second equation, and
then between the second and the third. This gives two equations
relating $U$ and $F_1$. Ignoring the irrelevant factors $U$, we then
eliminate $U$. This gives  
an algebraic equation satisfied by $F_1$,  containing three distinct
factors. The right one is easily identified, given that $F_1=1+O(t)$,
and one concludes that the \gf\ of quadrangular dissections of a
square, counted by the number of vertices, satisfies
$$F_1= 1-8t 
+2t \left(5- 6\,{t}\right) { F_1}
-2\,{t}^{2} \left( 1+3\,t \right) {{ F_1}}^{2}
-{t}^{4}{{ F_1}}^{3}.
$$

%%%%%%%%%%%%%%%%%%%%%%%%%%%%%%%%%%%%%%%%%%%%%%%%%
\section{A generic algebraicity theorem}
\label{section-generic}
%%%%%%%%%%%%%%%%%%%%%%%%%%%%%%%%%%%%%%%%%%%%%%%%%
%
Let $Q(y_0, y_1, \ldots , y_k, t, v)$ be a polynomial in $k+3$
indeterminates, with coefficients in a field $\GK$.
%  of characteristic $0$.
%, containing $\cs$. 
% In real life examples, this field is $\cs$, or
% $\cs(v_1, \ldots , v_m)$ for some indeterminates $v_i$. 
We consider the functional equation
\beq
\label{main-eq}
F(u)\equiv F(t,u) = F_0(u)+ t\  Q\Big( F(u), \Delta F(u), \Delta
^{(2)}F(u),\ldots ,\Delta ^{(k)}F(u), t, u\Big) ,
\eeq
where $F_0(u)\in \GK[u]$ is given explicitly  and the operator
$\Delta$ is the divided difference (or discrete derivative):
$$
\Delta F(u) = \frac{F(u)-F(0)}u.
$$
Note that 
$$
\lim _{u\rightarrow 0} \Delta F(u) =F'(0),
$$
where the derivative is taken with respect to $u$. The operator
$\Delta ^{(i)}$ is obtained by applying $i$ times $\Delta$, so that:
$$
\Delta ^{(i)} F(u) = \frac{F(u)-F(0)-uF'(0) -\cdots -
  u^{i-1}/(i-1)!\,F^{(i-1)}(0)}{u^i}. 
$$
Observe that all the equations met in Sections~\ref{section-intro}
to~\ref{section-examples} are of the 
form~\Ref{main-eq}, or can be easily transformed into an equation of
this form. Clearly,~\Ref{main-eq} has a unique
solution $F(t,u)$ in 
$\GK[u][[t]]$
% , the ring of formal power series in $t$ with polynomial
% coefficients in $u$ 
 (think of extracting from~\Ref{main-eq} the coefficient of $t^n$,
for $n=0,1, 2\ldots$). 
Upon multiplying~\Ref{main-eq} by a large power of $u$, one obtains a
polynomial equation of the form
$$
P\Big( F(u), F_1, \ldots , F_k, t,u\Big) =0,
$$
where   $F_i= F^{(i-1)}(0)/(i-1)!$ is the coefficient of $u^{i-1}$ in
$F(u)$, for $1\le i \le k$. 
Here is the main result of this section.
\begin{Theorem} 
\label{generic-thm}
The \fps \ $F(t,u)$ defined by~{\em \Ref{main-eq}} is algebraic over
$\GK(t,u)$.
\end{Theorem}
The proof requires the following  result~\cite[Prop.~X.8]{lang}.
\begin{Theorem}
\label{lang}
Let $\GK \subset \GL$ be a field extension. For $1\le i \le n$,  
let $P_i(x_1, \ldots ,
x_n)$ be a polynomial in $n$ indeterminates $x_1, \ldots , x_n$, with
coefficients in the (small) field $\GK$.  Assume  $F_1, \ldots , F_n$ are
$n$ elements of the (big) field $\GL$ that satisfy
$P_i(F_1, \ldots , F_n)=0$ for all $i\le n$.  
Let $J$ be the Jacobian matrix
$$J= \left( \frac{\partial P_i}{\partial x_j}(F_1, \ldots , F_n)\right)_{1\le
  i,j \le n}.
$$
If $\det(J) \not = 0$, then each $F_j$ is algebraic over $\GK$.
\end{Theorem}
\noindent
{\bf Proof of Theorem~\ref{generic-thm}.} The idea is of course to
apply the general strategy of Section~\ref{section-key}. However, in
order to avoid multiplicities in the roots $U_i$, we first introduce a
small perturbation of ~\Ref{main-eq}. Let $\eps$ be a
new indeterminate, and  consider the equation
\beq
\label{epsilon-eq}
G(u)\equiv G(z,u, \epsilon) = F_0(u)+ \eps ^k z \Delta ^{(k)}G(u) +
z^2 Q\Big( G(u), \Delta G(u), \Delta ^{(2)}G(u),\ldots 
,\Delta ^{(k)}G(u), z^2, u\Big) 
\eeq
where  $F_0$ and $Q$ are the same polynomials as above. Again, this
equation admits a unique solution in the ring of \fps \ in $z$ with
coefficients in $\GK[u,\eps]$. Moreover, $G(z,u, 0)= F(z^2,u)$, so
that it suffices to prove that $G(z,u, \epsilon)$ is algebraic over
$\GK(z,u,\eps )$.

\smallskip
We now apply to~\eqref{epsilon-eq} our  general strategy.
Our first task will be to convert~\Ref{epsilon-eq} into a polynomial
equation of the 
form~\Ref{eq-generic1}. Let $x_0, x_1, \ldots , x_k$ and $v$ be some
indeterminates. For $0\le i \le k$, let  
$$
Y_i= \frac{x_0-x_1-vx_2-\cdots - v^{i-1} x_i}{v^i}
$$
and let 
\beq \label{R-def}
R(x_0, x_1, \ldots , x_k, z, v)= x_0-F_0(v)- \eps ^k z Y_k
-z^2 Q\Big( Y_0, Y_1, \ldots , Y_k, z^2, v\Big)
 .
\eeq
Then 
$$
R\Big( G(u), G_1, \ldots , G_k, z,u\Big) =0,
$$
with $G_i= G^{(i-1)}(0)/(i-1)!$. Moreover, $R$ is a polynomial in $z$
and the $x_i$, but a rational function  in $v$.
So let $m$ be the smallest integer such that
\beq \label{P-R}
P(x_0, x_1, \ldots , x_k, z, v):= v^m R(x_0, x_1, \ldots , x_k, z, v)
\eeq
is a polynomial in $z, v$ and the $x_i$ (with coefficients in
$\GK(\epsilon)$). Then $m\ge k$ (because of the 
term $\eps ^k z Y_k$ occurring in $R$) and Eq.~\Ref{epsilon-eq} now reads
\beq
\label{eq-G-pol}
P\Big( G(u), G_1, \ldots , G_k, z,u\Big) =0.
\eeq
Let us apply to~\Ref{eq-G-pol} the general strategy of
Section~\ref{section-key}. We 
need to find sufficiently many fractional power series $U$ in $z$, with 
coefficients in some algebraic closure of $\GK(\eps)$, satisfying
$$
\frac{\partial P}{\partial x_0}\Big( G(U),  G_1, \ldots , G_k,
z,U\Big) =0.
$$
Let us focus on the non-zero solutions $U$. The above condition
is then  equivalent to
$$
U^k= \eps ^k z + z^2 
\sum_{i=0}^k U^{k-i} \, \frac{\partial Q}{\partial y_i}\Big(F(U), \ldots,
\Delta^{(k)} F(U), z^2, U\Big).
$$
By Theorem~\ref{thm-roots}, this equation has exactly $k$ solutions
$U_1, \ldots , U_k$, which are fractional power series in $z$ with
coefficients in an algebraic closure of $\GK(\eps)$. More precisely, the
Newton-Puiseux algorithm shows that these series  can be written as
\beq \label{Ui-expansion}
U_i= \eps \, \xi ^i  s\, (1+V( \xi ^i  s))
 \eeq
where $s=z^{1/k}$, $\xi$ is a primitive $k$th root of unity and $V(s)$
is a \fps \ 
in $s$ with coefficients in $\GK( \eps)$, having constant term $0$. In
particular, the $k$ series $U_i$ are distinct.

The following system of $3k$ polynomial equations thus holds:
$$
\forall i \in [1,k], \hskip 8mm \left\{ 
\begin{array}{lll}
P\Big(G(U_i), G_1, \ldots , G_k, z, U_i\Big) & = & 0,\\
%\frac{\partial P}{\partial x_0}
P'_0
\Big(G(U_i), G_1, \ldots , G_k, z, U_i\Big) & = & 0,\\
%\frac{\partial P}{\partial u}
P'_v
\Big(G(U_i), G_1, \ldots , G_k, z, U_i\Big) & = & 0,
\end{array} \right.
$$
where $P'_0$ and $P'_v$ respectively denote the derivatives of $P$
with respect to $x_0$ and $v$. 
The above system relates $3k$ unknowns, namely the $U_i$, the $G(U_i)$,
and the series $G_1, \ldots , G_k$, and has coefficients in $\GK(\eps,
z)$. Let us now apply Theorem~\ref{lang}. The Jacobian matrix is
represented below for $k=3$.  The rows are indexed by the $3k$
equations, and the columns by the $3k$ unknowns, taken in the following
order: $G(U_1), U_1, \ldots, G(U_k), U_k$ and finally $G_1, \ldots ,
G_k$.  We denote any series
of the form $S\left(G(U_i), G_1, \ldots , G_k, z, 
U_i\right)$ by $S(U_i)$ for short.
The notation $P'_i$ means that the derivative of $P$ is taken with
respect to the variable $x_i$.
$$
\left( 
\begin{array}{cc|cc|cc|ccccc}
%\displaystyle
%
%   Premiere ligne
%
%
P'_0 (U_1) & P'_v  (U_1)
& 0 & 0&0 &0 & P'_1 (U_1) & \cdots  & P'_k (U_1) \\
P''_{0,0} (U_1)  & P''_{0,v}(U_1)  & 0&0 &0 &0 &
\star & \cdots & \star\\
P''_{0,v}(U_1)  &
P''_{v,v}(U_1)  &0 &0 & 0& 0&
\star & \cdots & \star\\ 
&&&&&&&&\\ \hline
%
%
%   Deuxieme ligne
%
%
&&&&&&&&\\
0&0&P'_0 (U_2) & P'_v (U_2)
&0&0& P'_1(U_2) & \cdots  &
P'_k (U_2)  &\\
0&0&P''_{0,0}  (U_2)  & P''_{0,v}(U_2) &0&0&\star &
\cdots & \star\\ 
0&0&P''_{0,v}(U_2)  &
P''_{v,v}(U_2)&0&0&\star & \cdots & \star\\
&&&&&&&&\\\hline
&&&&&&&&\\
%
% Troisieme ligne
%
%
0&0 &0 &0&P'_0 (U_k) & P'_v(U_k)
  & P'_1(U_k) & \cdots  & P'_k (U_k) \\
0 &0 &0 &0&P''_{0,0} (U_k)  &
 P''_{0,v}(U_k)    &
\star & \cdots & \star\\
0&0 &0 &0& P''_{0,v}(U_k)  & P''_{v,v}(U_k)    &
\star & \cdots & \star\\ 
\end{array}
\right).
$$
Recall that 
\beq\label{U_j-char}
%\frac{\partial P}{\partial x_0} 
P'_0(U_j) = 
%\frac{\partial P}{\partial v} 
P'_v(U_j) =0
\eeq
for all $j$, so that the top line in each $3 \times 2$ rectangle is
actually zero. Consequently, the determinant factors into $k$ blocks
of size $2$ and one block of size $k$:
\beq
\label{factored-jacobian}
\det(J) = \pm \prod_{j=1}^k 
\left(   
%\frac{\partial ^2 P}{\partial x_0^2 }
P''_{0,0}( U_j) 
% \frac{\partial ^2 P}{\partial u^2}
P''_{v,v}(U_j)-
% \frac{\partial ^2 P}{\partial x_0 \partial u}
P''_{0,v}( U_j)^2\right)
\ \det\left(  
%\frac{\partial  P}{\partial x_j  }
P'_i ( U_j)\right)_{1\le i, j \le k}.
\eeq
Our aim is to prove that this Jacobian is not zero.

\medskip

\noindent 1. 
% Let us first prove that the term
Assume 
\beq\label{assumption}
P''_{0,0}( U_j) 
% \frac{\partial ^2 P}{\partial u^2}
P''_{v,v}(U_j)-
% \frac{\partial ^2 P}{\partial x_0 \partial u}
P''_{0,v}( U_j)^2 =0.
\eeq
% is not zero.
%
Let us differentiate twice the functional equation~\Ref{eq-G-pol} with
respect to $u$. We first obtain
$$
G'(u) P'_0(G(u), \ldots, u) + P'_v(G(u), \ldots, u)=0
$$
and then
$$
G''(u) P'_0( u) + G'(u)^2 P''_{0,0}( u) 
+ 2G'(u) P''_{0,v}( u)+ P''_{v,v}( u)=0,
$$
where, as above, the notation $S(u)$ actually stands for $S(G(u), G_1,
\ldots, G_k,t,u)$. For $u=U_j$, in view of~\Ref{U_j-char}, the latter
equation becomes 
$$
 G'(U_j)^2 P''_{0,0}(U_j) + 2G'(U_j) P''_{0,v}(U_j)
+ P''_{v,v}(U_j)=0.
$$
 The assumption~\Ref{assumption} implies that the
quadratic equation in $x$
$$
 x^2 P''_{0,0}(U_j) + 2x P''_{0,v}(U_j)+ P''_{v,v}(U_j)=0
$$
has a double root. The previous equation shows that this root is $
G'(U_j)$, so that
$$
G'(U_j)  P''_{0,0}(U_j)+ P''_{0,v}(U_j)=0.
$$
Given that $P'_0(U_j)=0$, this is equivalent to saying that the series 
$$
P'_0(G(u), G_1, \ldots, G_k, t, u)
$$
admits $u=U_j$ as a \emm multiple, root, whereas we have seen that the $k$
non-zero roots of this equation are distinct. We have thus obtained a
contradiction, and so~\Ref{assumption} cannot hold.
\medskip

\noindent 
2. Let us now focus on the second part of the
   expression~\Ref{factored-jacobian} of the Jacobian.
From~\Ref{P-R} and~\Ref{R-def}, we derive that for $j\ge 1$, and
   indeterminates $x_0, x_1, \ldots , x_k, z$ and $v$:
$$
P'_j(x_0, \ldots, x_k, z, v)= v^m 
R'_j(x_0, \ldots, x_k, z, v)
= -v^m \left( \eps^k z \,\frac{\partial Y_k}{\partial x_j}
+ z^2 \sum_{\ell =j}^k  \frac{\partial Y_\ell}{\partial x_j}
\,
Q'_\ell (Y_0, \ldots, Y_k,z^2,v)\right)
$$
where $Q'_\ell$ denotes the derivative of $Q(y_0, \ldots, y_k,t,v)$
with respect to  $y_\ell$. Given that
$$
 \frac{\partial Y_\ell}{\partial x_j} = -v^{j-\ell-1},
$$
the above derivative can be rewritten
\beq
\label{Ppj}
P'_j(x_0, \ldots, x_k, z, v)= v^{m-k} 
\left( \eps^k z v^{j-1}
+ z^2 \sum_{\ell =j}^k  v^{k-\ell + j-1}
Q'_\ell (Y_0, \ldots, Y_k,z^2,v) \right).
\eeq
Let us specialize this to 
$
P'_j(U_i) \equiv 
P'_j (G(U_i), G_1, \ldots , G_k, z, U_i)
$.
By~\Ref{Ui-expansion}, this is a \fps \ in $s=z^{1/k}$, with coefficients in
$\GK(\eps, \xi)$. Moreover, $z= s^k = o(U_i^{j-1})$
for $1\le j \le k$, so that, in view of~\eqref{Ppj}, the  first 
%non-zero 
term in the expansion of $P'_j(U_i)$ in $s$ is
$$
(\xi^i \eps s)^{m+j-1}.
$$
(Recall that $\xi^{ik}=1$.)
 The last factor in the determinant~\Ref{factored-jacobian} of the
Jacobian matrix $J$ reads 
$$
\det\left(  P'_i ( U_j)\right)_{1\le i, j \le k}=
\det\Big( (\xi^i \eps s)^{m+j-1} \Big)_{1\le i,j\le k} +
\hbox{higher powers of $s$}.
$$
But
$$
\det\Big( (\xi^i \eps s)^{m+j-1} \Big)_{i,j} = \prod_{j=1}^k (\eps s)^{m+j-1}
\prod_{i=1}^k (\xi^i)^{m} \det \Big( (\xi^i) ^{j-1}\Big)_{i,j}.
$$
The last term is the VanderMonde of the $\xi^i$. It equals
$$
\pm \prod _{1\le i<j\le k} \left( \xi^i -\xi ^j\right)
$$
and it is not zero, since $\xi$ is a $k$th primitive root of unity.

%Phew! 
We have at last proved that the determinant of the Jacobian matrix
associated with our system of $3k$ polynomial equations is not zero. By
Theorem~\ref{lang}, the series $G_i$ are algebraic over
$\GK(z,\eps)$. Recall that $G_i$ is, up to a multiplicative constant,
the derivative $G^{(i-1)}(0)$ of $G(u)$. In view of~\Ref{epsilon-eq},
the series $G(z,u,\eps)$ is 
algebraic over $\GK(z,u,\eps)$. By specializing $\eps$ to $0$, we
conclude that $F(t,u)=G(\sqrt t, u, 0)$ is algebraic over $\GK(t,u)$.
\cqfd

%%%%%%%%%%%%%%%%%%%%%%%%%%%%%%%%%%%%%%%%%%%%%%%%%
\section{Algebraicity results for planar maps}
\label{section-equations}
%%%%%%%%%%%%%%%%%%%%%%%%%%%%%%%%%%%%%%%%%%%%%%%%%

\begin{figure}[hbt]
\begin{center}
\epsfbox{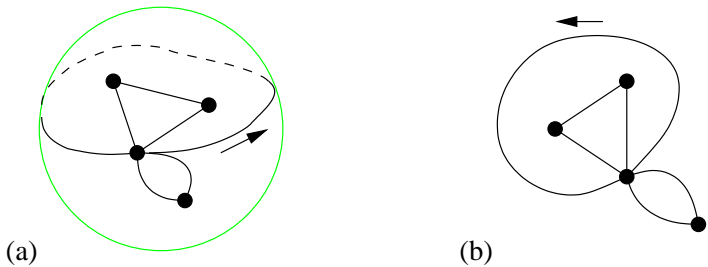}
\end{center}
\caption{(a) A rooted planar map on the sphere -- (b) Canonical
  representation on the plane.}
\label{planarmaps}
\end{figure}

A {\em planar map} is a 2-cell decomposition of the oriented
sphere into vertices (0-cells), edges (1-cells), and faces
(2-cells). Loops and multiple edges are allowed
(Figure~\ref{planarmaps}(a)).  The 
degree of a vertex (or a face) is the number 
of incidences of edges to this vertex (or face). 
Two maps are isomorphic if there exists an orientation
preserving homeomorphism of the sphere that sends cells of one of the
maps onto cells of the same type of the other map and preserves
incidences. We shall consider maps up to isomorphisms.

A map is {\em rooted} if one of its edges, called the root
edge,  is distinguished and oriented. In this case, the map can be
drawn in a canonical way in the plane, by deciding that the infinite
face lies to the right of the root-edge. This face is sometimes called
the root-face. Its degree is called the outer-degree. The 
starting point of the root-edge is  the  root-vertex.
A \emm  corner, of a face $F$ is a 3-tuple $(e_1,v,e_2)$, where $e_1$ and
$e_2$ are edges, $v$ is a vertex, and $e_1, v$ and $e_2$ are met
consecutively when 
walking around the face $F$ in counterclockwise order. The number of
corners of $F$ is thus its degree. In the map of
Figure~\ref{planarmaps}, the root-face has three corners.
In what follows, we  consider only  rooted maps,
and the word ``rooted'' is often omitted. 

A map $M$ is separable if it contains a vertex whose deletion
disconnects $M$. For instance, the map of Figure~\ref{planarmaps} is
separable, since deleting the root-vertex disconnects it.

\begin{figure}[hbt]
\begin{center}
\epsfbox{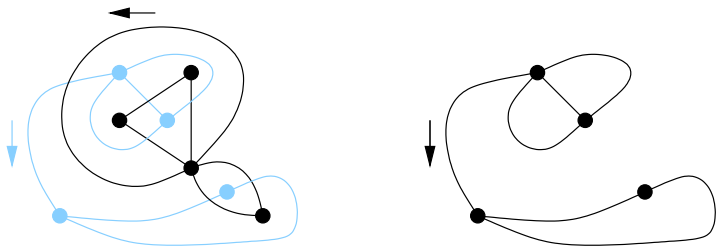}
\end{center}
\caption{Construction of the dual map.}
\label{fig-dual}
\end{figure}

 The dual map $M^*$ of a map $M$ describes the 
incidence relation between the faces of $M$ (Figure~\ref{fig-dual}). To
  construct $M^*$, create 
  a vertex  in every face of $M$: this gives the vertices of
  $M^*$. The edges of $M^*$ are in bijection with the edges of $M$:
  for each edge $e$ of $M$, incident to the faces $f_1$ and $f_2$, create
  an edge of $M^*$ that crosses $e$ and joins the vertices of $M^*$
  corresponding to $f_1$ and $f_2$. The root-edge of $M^*$ is chosen
  canonically. 

%%%%%%%%%%%%%%%%%%%%%%%%%%%%%%%%%%%%%%%%%%%%%%%%%
\subsection{The face-distribution of planar maps}
\label{subsection-distribution}
\begin{figure}[hbt]
\begin{center}
\input{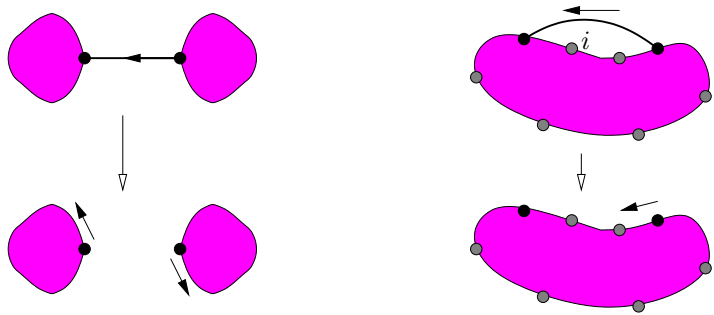}
\end{center}
\caption{The decomposition of planar maps.}
\label{fig:face-distrib}
%\hrule
\end{figure}

%%%%%%%%%%%%%%%%%%%%%%%%%%%%%%%%%%%%%%%%%%%%%%%%%
Many functional equations for
planar maps are based on the deletion of the root-edge. Here, we write
an equation for the series $F(t,u;z_1, \ldots , z_m, 
\ldots)=F(t,u;\bm z)$ that counts rooted planar maps by the number of edges
(variable $t$), the outer-degree (variable $u$) and the
number of \emm finite, faces of degree $i$ (variable $z_i$) for all $i
\ge 1$. This equation essentially appears in an old paper of
Tutte~\cite[Eq.~(1)]{tutte-general}.  
\begin{Lemma}
\label{lemma-finite-faces} 
The \gf\ $F(t,u;\bm z)\equiv F(u)$
satisfies
$$
F(u)=1+tu^2F(u)^2+t \sum_{i\ge 1} z_i\, \frac{F(u)-\sum_{j=0}^{i-2} u^j
  F_j}{u^{i-2}} ,
$$
where $F_j$ is the coefficient of $u^j$ in $F(u)$. 
\end{Lemma}
\noindent
{\bf Proof.} Take a planar map $M$. If it is not reduced to a single
vertex, delete the root-edge (but not its endpoints). Then
\begin{itemize}
\item[--] either two connected components are left, which we can root
  in a canonical way (Figure~\ref{fig:face-distrib}). The \gf\ of
  such maps is $tu^2F(u)^2$,
\item[--] or only one connected component is left, which we can root
  in a canonical way. Let $j$ be its outer-degree, and let $i$ be the
  degree of the finite face that has been deleted with the root-edge of
  $M$. Then $i\in[1, j+1]$. The \gf\ of maps of this second type is 
$$
t\sum_{j\ge 0}\left( F_j \sum_{i=1}^{j+1} z_i u^{j-i+2}\right).
$$
\end{itemize}
Adding the two contributions gives a functional equation for $F(u)$ which 
\begin{itemize}
\item[--] specializes to~{\rm \Ref{eq-planar-maps}}  when $z_i=1$ for all $i$,
\item[--]  gives the equation of Lemma~\ref{lemma-finite-faces} upon
exchanging the order of the summations on $i$ and $j$.
\end{itemize}
\cqfd

One may think that there is in $F(t,u; \bm z)$ an unpleasant lack of
symmetry: why should one count only the \emm finite, faces of a given
degree? Let  $G(t;z_1, \ldots , z_m,
\ldots)=G(t;\bm z)$ count rooted planar maps by the number of edges
(variable $t$) and the
number of faces (finite or not) of degree $i$ (variable
$z_i$). Observe that, by duality, $G(t;
\bm z)$  also counts planar 
maps by the number of edges and the number of \emm vertices, of degree
$i$. We call $G$ the \emm face-distribution \gf, of planar maps (equivalently, the
vertex-distribution \gf\ of planar maps).
\begin{Lemma}
\label{lemma-all-faces} The face-distribution \gf\ of planar maps,
$G(t;\bm z)$, is related to the series  $F(t,u;\bm z)$ of
Lemma~{\rm \ref{lemma-finite-faces}}  by 
$$
G(t; \bm z)=\frac 1 t [u^2] F(t,u;\bm z).
$$
\end{Lemma}
\noindent
{\bf Proof.} Take a map $M$ with outer-degree 2. The root-face is
incident to two edges: delete the non-root one to obtain a planar
map $M'$. This transformation is bijective and the degree distribution of
\emm finite, faces in $M$ coincides with the degree distribution of
all faces in $M'$.
\cqfd
The equation of Lemma~\ref{lemma-finite-faces} was
 solved in~\cite{bender-canfield} in the case where $z_i=1$ if $i\in D$ and
$z_i=0$ otherwise, for a given set $D$. More recently, the
 vertex-distribution \gf\ of planar maps was characterized  
 in~\cite{BDG-planaires} 
via two methods: first, by a matrix integral calculation,  and then
using a purely bijective approach. In
 Section~\ref{section-distribution}, we  provide an alternative
solution to this problem, and prove that it is
equivalent to~\cite{BDG-planaires}.
   For the moment, observe that the generic algebraicity
theorem of Section~\ref{section-generic}
(Theorem~\ref{generic-thm}) implies the following:
\begin{Coro}\label{coro-face-algebraic}
Let $m\ge 1$, and let $F(t,u;z_1, \ldots , z_m)$ be the \gf\ of rooted
planar maps in which no \emm finite, face has a degree larger than $m$ (as
above, $t$ counts edges, $u$ the outer-degree, and $z_i$ the number of
finite faces of degree $i$). Similarly, let  $G(t;z_1, \ldots , z_m)$
be the face-distribution \gf\ of rooted 
planar maps in which no face has a degree larger than $m$. Then both
series are algebraic.
\end{Coro}
\noindent{\bf Proof.} These series $F$ and $G$ are obtained by setting
$z_i=0$ for all $i>m$ in the series $F$ and $G$ of
Lemmas~\ref{lemma-finite-faces} 
and~\ref{lemma-all-faces}. The equation of Lemma~\ref{lemma-finite-faces} has then the generic
form~\Ref{main-eq}. By Theorem~\ref{generic-thm}, its solution
$F(t,u;z_1, \ldots , z_m)$ is algebraic over $\qs(t,u,z_1, \ldots ,
  z_m)$. Since the extraction of coefficients preserves algebraicity,
Lemma~\ref{lemma-all-faces} implies that $G(t;z_1, \ldots , z_m)$ is
algebraic too. 
\cqfd

%%%%%%%%%%%%%%%%%%%%%%%%%%%%%%%%%%%%%%%%%%%%%%%%%
\subsection{The face-distribution of Eulerian planar maps}
\label{section-eulerian}
%%%%%%%%%%%%%%%%%%%%%%%%%%%%%%%%%%%%%%%%%%%%%%%%%
The question we address here is similar to that of
Section~\ref{subsection-distribution}, but 
is made harder by the fact that we now deal with \emm Eulerian, maps, that
is, with maps in which all vertices have an even degree.
The faces of an Eulerian map can be uniquely coloured in black and
white in such a way 
\begin{itemize}
\item[--] the infinite face is white,
\item[--] every black face is only adjacent to white faces, and
  vice-versa.
\end{itemize}
Let $F(t,u; x_1, x_2, \ldots ; y_1, y_2, \ldots) = F(t,u;\bm x, \bm y)$
 be the \gf\ of these maps, where $t$ counts edges, $u$ the
outer-degree, $x_i$ the number of (finite) white faces of degree $i$,
and $y_i$ the number of black faces of degree $i$ (all
black faces are finite).

If we set $y_i=0$ for $i\not = 2$, the series $ F(t,u;\bm x, \bm y)$
only count those Eulerian maps in which every black face has degree
2. Contracting every black face into a single edge gives a planar map
whose face-distribution coincides with the white face-distribution of
the original Eulerian map. Consequently,  $F(t,u; z_1, z_2,
\ldots ; 0, 1, 0, \ldots)$ is the series studied in
Lemma~\ref{lemma-finite-faces}, and 
the problem addressed here  generalizes the previous one.
% of the face-distribution enumeration of (general) planar maps.

%Similarly, if we set $y_i=0$ for $i\not =m$, for some given $m\ge 2$,
%and $x_i=0$ when $i$ is not a multiple of $m$, we are back to the
%enumeration of $m$-constellations (Section~\ref).

\begin{figure}[hbt]
\begin{center}
\input{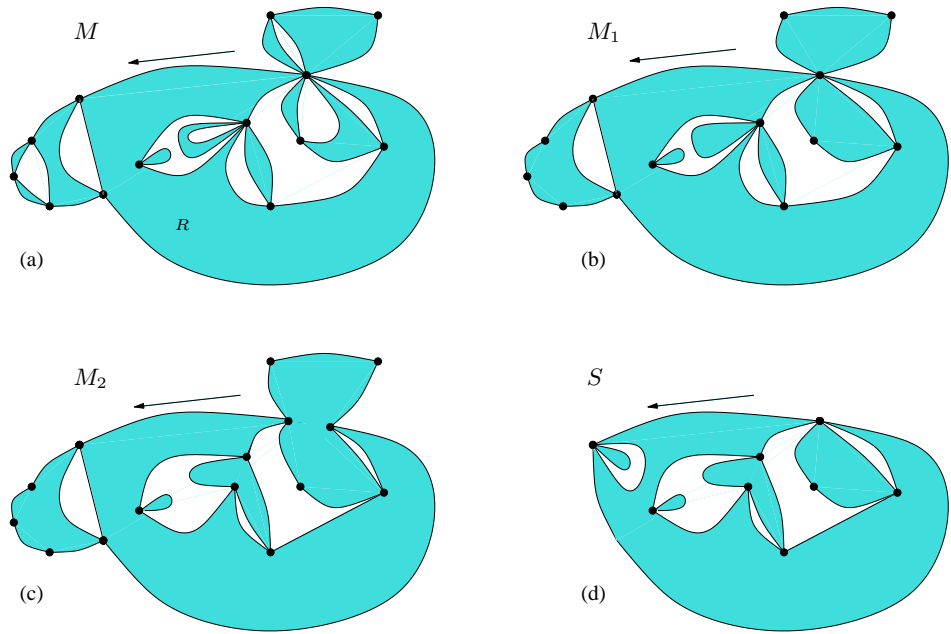}
\end{center}
\caption{From an Eulerian map $M$ to a skeleton $S$.  In  step (b), all
  white faces that are not adjacent to the black root-face $R$
  disappear. In  step (c), all connected component that are left
  after deleting the black root-face are inflated to
  polygons. Finally, in step (d), all edges that are incident to the
  white root-face but not to the black root-face are contracted.}
\label{fig:eulerian}
%\hrule
\end{figure}

In order to obtain a functional equation for $F(t,u;\bm x, \bm y)$, we
will delete all the edges of the black face
incident to the root-edge. 
We call this face the \emm black
root-face,. A face is called a \emm polygon, if the number of vertices
it contains coincides with its degree. 
\begin{Definition}
An Eulerian map $M$ is a \emm skeleton, if the following conditions
hold:
\begin{itemize}
\item[$(i)$] each of the connected
components that remain after deleting the edges of the black root-face $R$
is  either a single vertex or a polygon,
\item[$(ii)$] every edge that is incident to the white root-face is also
  incident to the black root-face.
\end{itemize}
A connected component of $M\setminus R$  is
called an \emm internal component of $M$, if none of its vertices
belong to the infinite  
face. Otherwise, it is said to be an \emm external component of $M$.,
\end{Definition}
The fourth map of Figure~\ref{fig:eulerian} is a skeleton. Among its
non-root black faces, two are external, and two are internal.
The following observation will be useful to prove that the
face-distribution \gf\ of Eulerian maps with faces of bounded degree
is algebraic.
\begin{Lemma}
\label{lemma-finite-skeletons}
Let $m\ge 1$. There exists only a finite number of skeletons in which
the black root-face and all the
finite white faces have degree at most $m$.
\end{Lemma}
\noindent{\bf Proof.}
Let us first bound the number of white faces. Condition $(i)$ implies
that  each  white face of a skeleton shares at
least one edge with the black root-face.  Conversely, each edge of
the black root-face belongs to exactly one white
face. Since there are, by assumption, at most $m$ such edges, the
number of white faces is at most $m$. 
By assumption, the finite white faces have degree at most
$m$. Condition $(ii)$ implies that this is also true for the infinite
white face. Consequently, the total number of
edges that are incident to a  white face --- that is, the total number
of edges --- is at most $m^2$.
Since there only exists a finite number of 
maps having a given number of edges,  the result follows.
\cqfd
\begin{Propo}
\label{propo-finite-faces-eulerian} 
Let  ${\mathcal S}$ denote the set of skeletons. The \gf\ $F(t,u;\bm
x;\bm y)\equiv F(u)$ counting Eulerian maps 
according to the above-defined parameters satisfies
$$
F(u)=1+ 
\sum_{S \in {\mathcal S}}
\left( u^{d(S)}  t^{i(S)}y_{i(S)} \prod_{k\ge 1} x_k^{w_k(S)}
\prod_{k\ge 1} F_k^{I_k(S)} 
\prod_{k\ge 0}
%  \left( \frac{F(u)-\sum_{j=0}^{k-1} u^j   F_j}{u^{k}}\right)
\left( \Delta^{(k)}F(u)\right)^{E_k(S)} \right),
$$
where   
% ${\mathcal S}$ denotes the set of skeletons,
% in which the black root-face has degree $i$, and
 for any skeleton $S$, $d(S)$ is the outer-degree, $i(S)$ is the
 degree of the black 
 root-face, 
$w_k(S)$ is the number of finite white faces of
degree $k$, and $I_k(S)$ (resp. $E_k(S)$) is the number
of internal 
(resp. external) components of degree $k$.
As above, $F_j$ denotes the coefficient of $u^j$ in $F(u)$, 
and for $k \ge 0$,
$$
\Delta^{(k)}F(u)=\frac{F(u)-\sum_{j=0}^{k-1} u^j   F_j}{u^{k}}.
$$
%This equation specializes to~\Ref{eq-planar-maps}  when $z_i=1$ for all $i$.
\end{Propo}
\noindent
{\bf Proof.} Take an Eulerian map $M$, not reduced to a single
vertex. We first describe how to
associate a skeleton to $M$. This construction is illustrated in
Figure~\ref{fig:eulerian}. 
 Let $R$ denote the black root-face of $M$. 
% Let $i$ be the degree of $R$. Label the corners of $R$ with the
% labels 1,2, \ldots , $i$, in counterclockwise order. Delete the
% edges of $R$, keeping the vertices  and the labels. 
Consider the set of connected components that are left after the
deletion of the edges of $R$ 
(since we do not delete the vertices of $R$, some of these
components may be reduced to a single vertex). 
The corresponding
sub-maps of $M$ are called, for short, the \emm  
components, of $M$. Each component is itself  an Eulerian map. In order
to obtain a skeleton, we are going to modify the components of $M$,
while keeping the black root-face $R$ unchanged. In each
component, delete every edge that is
not in the infinite face of $M\setminus R$: in the resulting map $M_1$,
every component has only black (finite)
faces (Figure~\ref{fig:eulerian}(b)). Then ``inflate'' each component
into a black polygon having the same 
outer-degree (Figure~\ref{fig:eulerian}(c)).
%Glue back the set of polygons thus obtained to the black root-face
%$R$, following the indications given by the labels. 
This gives an
Eulerian map $M_2$. Finally, contract every edge of $M_2$ that is incident to
the white root-face but not to the black root-face. This gives a
skeleton $S$ (Figure~\ref{fig:eulerian}(d)). The finite white faces
of $S$ are in one-to-one 
correspondence with the finite white faces of $M$ that are adjacent
to $R$, and this correspondence preserves the degree.

Conversely, take a skeleton $S$ with black root-face of degree $i$. We
wish to find the \gf \ of Eulerian maps $M$ associated with $S$. To obtain
these maps, one must:
\begin{itemize}
\item[--] Replace every internal component of degree $k$ by an
  Eulerian map of outer-degree $k$; this gives the factors $F_k$ in
  the functional equation of Proposition~\ref{propo-finite-faces-eulerian}.
\item[--]  Replace every external  component  of degree $k$ by an
  Eulerian map of outer-degree $j\ge k$. Then $j-k$ edges of this map
  contribute to the outer-degree of the final map $M$. This gives the
  factors $\Delta^{(k)}F(u)$ in the equation.
\end{itemize}
The remaining factors take care of  $R$ and its edges, and of the
contribution of the white faces of $S$. The result follows.
\cqfd

\begin{Coro}
Let $m\ge 2$. Let $F(t,u;x_1, \ldots , x_m; y_1, \ldots , y_m)$ be
the \gf \ of Eulerian planar maps whose \emm finite, faces have degree
at most $m$, 
counted, as above, by the number of edges, the outer-degree, and the
degree-distribution of  black and white \emm finite, faces.

Similarly, let  $G(t,u;x_1, \ldots , x_m; y_1, \ldots , y_m)$  be
the \gf \ of Eulerian planar maps in which all  faces have degree
at most $m$, 
counted  by the number of edges, the outer-degree, and the
degree-distribution of white and black  faces.

Then $F$ and $G$ are algebraic.
\end{Coro}
\noindent{\bf Proof.} The  series $F$ is obtained by setting
$x_i=y_i=0$ for all $i>m$ in the series  of
Proposition~\ref{propo-finite-faces-eulerian}. In the 
equation given in  this
proposition, it is clear that the skeletons in which either the black
root-face, or one of the finite white faces, has degree more than $m$,
have a zero contribution. By Lemma~\ref{lemma-finite-skeletons}, the
right-hand side of the 
functional equation contains only finitely many terms, so that one can
apply Theorem~\ref{generic-thm}, and conclude that $F(t,u;x_1,
\ldots , x_m; y_1, \ldots , y_m) $ is algebraic.

In particular, the coefficient of $u^i$ in this series is algebraic.
Given that
$$G(t,u;x_1, \ldots , x_m; y_1, \ldots , y_m)= \sum_{i=0}^m 
x_i u^i \,[u^i]F(t,u;x_1,\ldots , x_m; y_1, \ldots , y_m),
$$
the algebraicity of $G$ follows.
\cqfd

%\medskip
\newpage
\begin{floatingfigure}[r]{4cm}
\flushright{\epsfbox{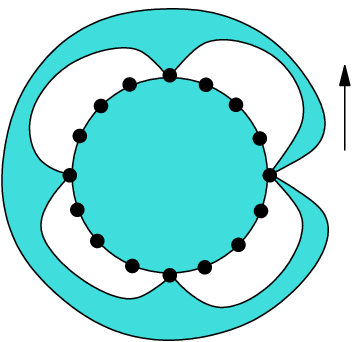}}
%\caption{gna}%{ skeleton}
%\label{fig-flower}
\end{floatingfigure}
\noindent{\bf Note.} It was already proved
 in~\cite{mbm-schaeffer-ising} that $F_2$, the
coefficient of $u^2$  in the series $F(t,u;\bm x, \bm y)$, is
algebraic.  The above corollary thus extends this earlier
result, and actually seems difficult to obtain via the combinatorial
 approach of~\cite{mbm-schaeffer-ising}. However, as far as $F_2$ is
 concerned, the result  of~\cite{mbm-schaeffer-ising} is more precise
 than a simple algebraicity statement,  since a system 
of  $2m+3$ polynomial equations defining $F_2$ is given explicitly,
together with its combinatorial interpretation. Let us compare
 the size of this system with the number of unknown series in our
 functional equation.  The skeleton of an
Eulerian map in which all finite faces have degree at most $m$ may
contain a component of degree $(m-1)^2$ (see
the figure for an example with $m=5$), but no more, so that the  functional
equation contains approximately $m^2$ unknown
functions. Consequently, the size of the polynomial system given by our 
general strategy is \emm quadratic, in $m$.

\begin{figure}[hbt]
\begin{center}
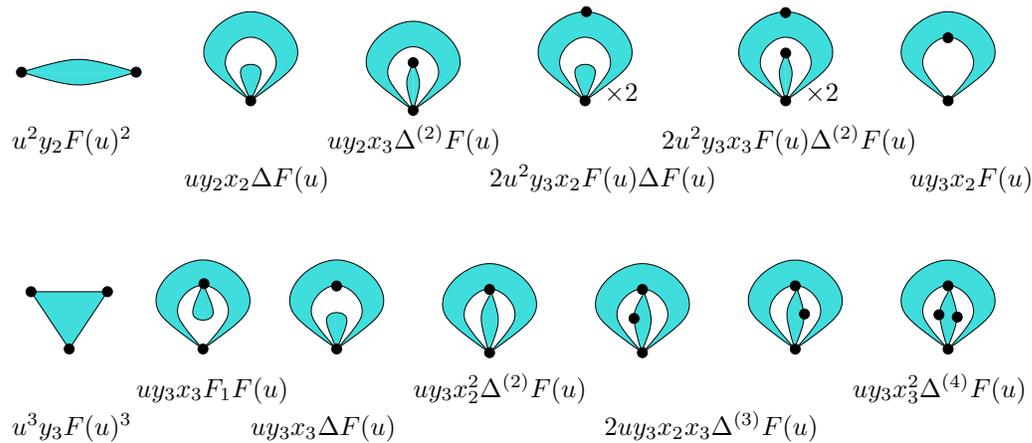
\end{center}
\caption{The skeletons that are involved in the enumeration of
  Eulerian maps with finite faces of degree 2 and 3. The multiplicities
  account for the number of possible rootings. }
\label{fig-23skeletons}
\end{figure}

\bigskip
\noindent{\bf Example.}
Let us illustrate Proposition~\ref{propo-finite-faces-eulerian} by
writing a functional equation 
for the \gf\ of Eulerian maps in which all finite faces have degree 2 or
3. The corresponding skeletons are shown in
Figure~\ref{fig-23skeletons}.
 Proposition~\ref{propo-finite-faces-eulerian} gives the
contribution  of each skeleton in the functional equation (for the sake
of simplicity, the variable $t$ is omitted: it is easily
recovered upon replacing $y_i$ by $t^iy_i$).

\noindent The functional equation reads:
$$
F(u)=1+ u^2y_2F(u)^2+u^3y_3F(u)^3
%
%+ux_2y_3F(u)+ux_3y_3F_1F(u)
%
+uy_3(x_2+x_3F_1)F(u)
+u\left(x_2y_2+2ux_2y_3F(u)+x_3y_3 \right) \Delta F(u)
$$
$$
+ u \left(x_3y_2+2ux_3y_3F(u)+x_2^2y_3 \right) \Delta^{(2)} F(u)
+2ux_2x_3y_3\Delta^{(3)} F(u)
+ ux_3^2y_3\Delta^{(4)} F(u).
$$
We may check the validity of this equation as follows. 
Replacing $y_i$ by $t^iy_i$, we derive from this equation the
first terms of the expansion in $t$ of $F(u)$. Retaining only the
coefficient of $u^2$, we obtain the expansion of the series $F_2$ that
counts maps of outer-degree 2, and we  check that this expansion
is (fortunately!) in adequation with the algebraic equations
of~\cite{mbm-schaeffer-ising}.

%%%%%%%%%%%%%%%%%%%%%%%%%%%%%%%%%%%%%%%%%%%%%%%%%
\subsection{Constellations}
\label{subsection-constellations}
\begin{figure}[hbt]
\begin{center}
\epsfbox{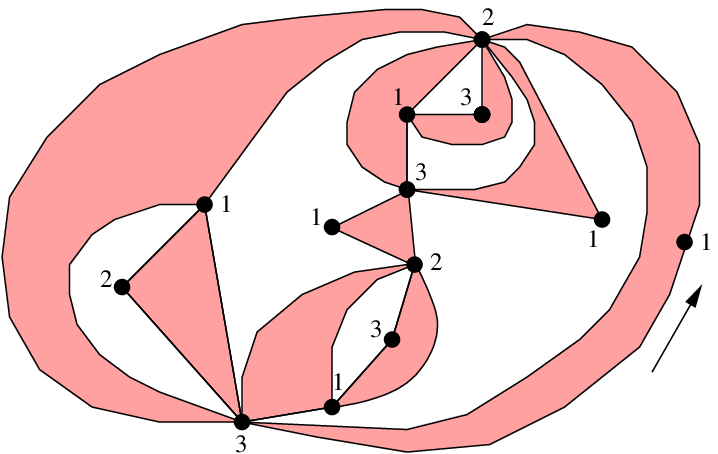}
\end{center}
\caption{A 3-constellation with its canonical labelling of root $3$.}
\label{fig-constellation}
\end{figure}

%%%%%%%%%%%%%%%%%%%%%%%%%%%%%%%%%%%%%%%%%%%%%%%%%
We focus in this section on the enumeration of certain Eulerian planar
maps defined by constraints on their face degrees. 
Let $m\ge2$.  An Eulerian planar map $M$, having its faces bicolored
in such a way the infinite face is white, is an \emm
$m$-constellation, if
\begin{itemize}
%\item[--] the infinite face is white,
%\item[--] every white face is only adjacent to black faces, and
%  vice-versa, 
\item[--] the degree of every black face is $m$,
\item[--] the degree of every white face is a multiple of $m$.
\end{itemize}
An example of a $3$-constellation is given in
Figure~\ref{fig-constellation}. As explained
in~\cite{mbm-schaeffer-constellations}, these maps 
are closely connected to \emm minimal transitive factorizations, of
permutations.

%The existence of a two-colouring of
%the faces is equivalent to the fact that every vertex has an even
%degree: these maps is \emm Eulerian,. Also, 
The above conditions
guarantee that it is possible to label all vertices, with labels taken
from the set $\{1,2,\ldots,m\}$, in such a way that in every black
face, the vertices are labelled $1,2, \ldots, m$ in counterclockwise
order. Moreover, if we fix the label of the root-vertex to be $i$, then there
is a unique labeling satisfying the above property, which we call the
\emm canonical labeling of root $i$., 

Let 
\beq
\label{GF-constellations-def}
F(t,u)\equiv F(u)= \sum_{n,d} a_{n,d} t^n u^d=\sum_{d} F_du^d,
\eeq
where $a_{n,d}$ is the number  of $m$-constellations having $n$ black
faces and  outer-degree $md$. This series is  a specialization of the
face-distribution \gf\ of Eulerian planar maps studied in
Section~\ref{section-eulerian}. More precisely, if, in the series
$F(t,u;\bm x , \bm y)$ of 
Proposition~\ref{propo-finite-faces-eulerian}, we set
$$
\left\{
\begin{array}{llll}
x_i&=& 1 &\hbox{if } m \hbox{ divides } i,\\
x_i &=& 0 &\hbox{otherwise},
\end{array}\right.
\quad \quad \hbox{and} \quad \quad 
\left\{\begin{array}{llll}
y_m&=& 1 ,\\
y_i &=& 0 & \hbox{if } i\not = m,
\end{array}\right.
$$
we obtain the series $F(t^m, u^m)$, with $F(t,u)$ defined
by~\Ref{GF-constellations-def}. However, the functional equation of
  Proposition~\ref{propo-finite-faces-eulerian}, specialized to the
  above values of $x_i$ and $y_i$, contains infinitely
  many terms.
We  give in Proposition~\ref{propo-constellations} an 
equation with \emm finitely, many terms defining $F(t,u)$.
 Before we do so, let us examine 
%the case $m=3$.
cases $m=2$ and $m=3$.

\smallskip
\noindent 
{\bf 2-Constellations.} Take a 2-constellation not reduced to a
single vertex, label the root-vertex with 2 and the other vertices
canonically. Each 
black face has degree 2 and contains a vertex labelled 1 and a vertex
labelled 2.  Contract each  black face
 to a single edge: this gives a \emm  bipartite, map,
that is, a map in which every face has an even degree. The series $F(t,u)$
thus counts bipartite maps by the number of edges ($t$) and half the
outer-degree (in other words, the number of corners labelled 1 in the
infinite face). Deleting the root-edge as we did in
Section~\ref{subsection-distribution} for general maps
% (Figure~\ref) 
now gives
\begin{eqnarray}
F(u) &=&1 + tuF(u)^2 + t \sum_{d\ge 0} F_d \left( u^d +\cdots +
u\right)\nonumber\\
&=& 1  + tuF(u)^2 +tu\, \frac {F(u)-F(1)}{u-1}.\label{eq-bipartite}
\end{eqnarray}
Observe that the deletion of the root-edge in a bipartite map
corresponds to the deletion of the black root-face
% adjacent to the root-edge
in the associated 2-constellation. The study of
2-constellations will be useful in
Section~\ref{section-hard-equation}, where we count 
certain maps with bicolored vertices. However, it is a bit  too simple
to foresee what 
 happens for general $m$-constellations. This is why we also treat below
the case of 3-constellations. 

\smallskip
%%%%%%%%%%%%%%%%%%%%%%%%%%%%%%%%%%%%%%%%%%%%%%%%%%%%%%%%%%%%%%%%%

\begin{figure}[thb]
\begin{center}
\input{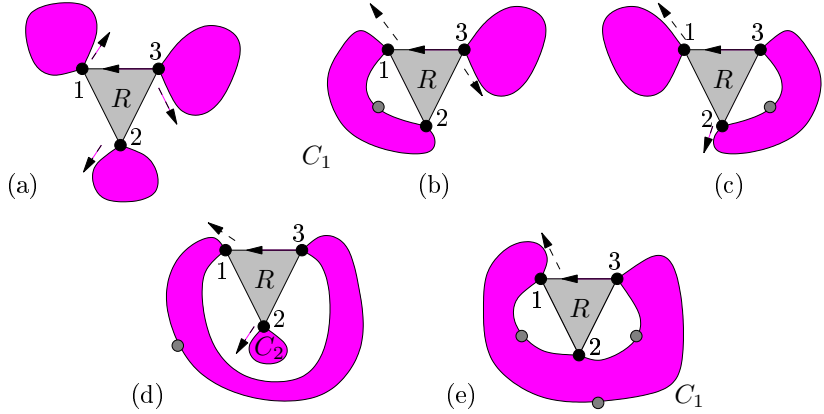}
\end{center}
\caption{The decomposition of 3-constellations. The dashed arrows
  indicate how to root the components after the deletion of the black
  root-face. }
\label{fig:3-const}
%\hrule
\end{figure}

\noindent 
{\bf 3-Constellations.} Take a 3-constellation $C$ not reduced to a
single vertex, label the root-vertex with 3 and the other vertices
canonically. Let $R$ denote the black root-face.
% that is adjacent to the root-edge. 
Erase all the edges of  $R$ (but not its
vertices).  This leaves a set of connected components, which
are constellations, and  which we root
in a canonical way (Figure~\ref{fig:3-const}).

Five cases occur, depending on which vertices of $R$ end up in the
same component. For the first case, the \gf\ is
clearly $tuF(u)^3$. The second and third cases are symmetric and thus
give the same \gf . Note that the component $C_1$
 in Figure~\ref{fig:3-const}(b) must have outer-degree 3 at least, and that the
 number of ways to glue a (rooted) 3-constellation $C_1$ of
 outer-degree $3d$ to the
 face $R$ is $d$. If the $j$th corner labelled 2 of the infinite face
 of $C_1$ is glued to $R$, then $1+3(j-1)$ edges of $C_1$ contribute
 to the outer-degree of $C$.  Thus the \gf\ in the second case is
$$
t u^{2/3} F(u) \sum_{d\ge 1}\left( F_d \sum_{j=1}^d u^{1/3 +j-1}
\right)= tu F(u)\, \frac {F(u)-F(1)}{u-1}.
$$
In the fourth case, the component $C_2$ does not contribute to the
outer-degree of $C$, but this case is otherwise similar to the
previous one.  The \gf\ is now
$$
 tu F(1)\,\frac {F(u)-F(1)}{u-1}.
$$
Finally, in the fifth case, the component $C_1$ has degree $3d$ with
$d\ge 2$. Assume the $j$th corner labelled 3 of the infinite face of
$C_1$ is glued to $R$, as well as the $k$th corner 
labelled 2. 
Then $1\le j <k \le d$ and the \gf\ of this last case is
$$
tu^{1/3} \sum_{d\ge 2} \left( F_d u^{2/3} \sum_{j=1}^{d-1}\left( u^{j-1}
\sum_{k=j+1}^{d} 1\right)\right),
$$
which, after two summations, reduces to
$$
tu \frac {F(u)-F(1)-(u-1)F'(1)}{(u-1)^2}.
$$ 
Finally, the \gf\ of 3-constellations satisfies
\beq
\label{eq-3const}
F(u)=1+tuF(u)^3+ tu (2F(u)+F(1)) \,\frac {F(u)-F(1)}{u-1} 
%+ tu F(1)\,\frac {F(u)-F(1)}{u-1}
+tu \frac {F(u)-F(1)-(u-1)F'(1)}{(u-1)^2}.
\eeq

\smallskip

\noindent 
{\bf $m$-Constellations.} In  order to write a functional equation for
general $m$-constellations, we need the notion 
of \emm non-crossing partitions,~\cite{rodica-noncrossing}. A partition $P$ of the set
$\{1,2,\ldots,m\}$  is \emm non-crossing, if one cannot find
$i<j<k<\ell$ such that $i$ and $k$ are in the same block, $j$ and $\ell$
are in the same block, but $i$ and $j$ are \emm not, in the same block.
%(Figure~\ref{fig-ncpartition}). 
A block $B$ of a non-crossing
partition $P$ is \emm internal, if
there exists another block $B'$ such that $\min B' < \min B \le \max B
<\max B'$. Otherwise, it is external. Let  ${\mathcal P}_m$ denote the set of
non-crossing partitions of
$\{1,2,\ldots,m\}$.
\begin{Propo}
\label{propo-constellations}
Let $m\ge 2$. The \gf\ $F(t,u)\equiv F(u)$ of $m$-constellations,
defined by~{\em\eqref{GF-constellations-def}}, satisfies:
$$
F(u)=1+ tu \sum_{P\in {\mathcal P}_m} 
\prod_{k=1}^{m-1}   \left( G_{k-1}\right)^{I_k(P)}
\prod_{k=1}^{m}
\left( \frac {F(u) -\sum_{i=0}^{k-2} {(u-1)^i}
G_i} {(u-1)^{k-1}}\right)^{E_k(P)},
$$
where 
$$
G_i=\frac 1 {i!}\, \frac{\partial ^{i} F}{\partial u^i}(1)
$$
and $I_k(P)$ (resp. $E_k(P)$) denotes the number of internal
(resp. external) blocks of cardinality $k$ in the partition $P$.
\end{Propo}
\noindent
Note that
$$
 G_{k-1}=\lim_{u\rightarrow 1}  \frac {F(u) -\sum_{i=0}^{k-2} {(u-1)^i}
G_i} {(u-1)^{k-1}}.
$$
The above equation defining $F(u)$ has degree $m$ in $F(u)$ and
involves $m-1$ additional unknowns series $G_i$, for $0\le i\le m-2$.\\
{\bf Proof.} The proof is based again on the deletion of the black
root-face. We  call this face the \emm root 
$m$-gon, and denote it by $R$. 

\begin{figure}[thb]
\begin{center}
\input{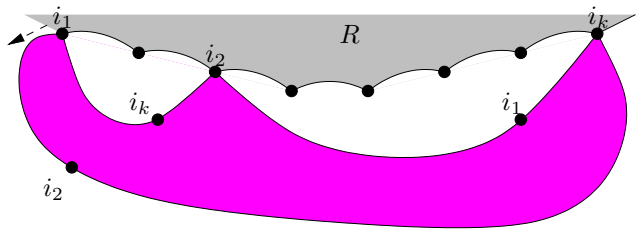}
\end{center}
\caption{The decomposition of $m$-constellations. The dashed arrow
  indicates how to root the component after the deletion of the black
  root-face. One has $i_1 <i_2< \cdots < i_k$.}
%The block $\{i_1, i_2, \ldots , i_k\}$ is associated to   the component}
\label{fig:decomp-const}
%\hrule
\end{figure}

\noindent {\bf 1. Decomposition of constellations.} Take a
constellation $C$ that is not reduced to a single vertex. Label the
root-vertex by $m$, and all the other vertices in a canonical
way. Erase all the edges of the root $m$-gon $R$ (but not its
vertices). This leaves a number of constellations, which we root
in a canonical way (Figure~\ref{fig:decomp-const}). For each of them, the
label of the 
root-vertex is minimal among the labels of the vertices that it
shares with $R$.

Associate with $C$ the partition $P$ of $\{1,2,\ldots,m\}$ defined as
follows:  $i$ and $j$ belong to the same block if and only if the
vertices labeled $i$ and $j$ in $R$ end up in the same connected
component after deleting the edges of $R$. By planarity of $C$, the
partition $P$ is non-crossing. To each block $B$ of $P$, there corresponds
a constellation $C_B$ (the associated connected component). 
% Its outer-degree is at least $m(|B|-1)$.

The outer-degree of $C$ is
\beq
\label{degree}
\hbox{ext}(P)+ \sum_{B \hbox{ external}} \delta(C_B,C),
\eeq
where $\hbox{ext}(P)$ is the number of external blocks of $P$ and
$\delta(C_B,C)$ is the number of edges of $C_B$ that contribute to the
outer-degree of $C$.  

\noindent {\bf 2. Construction of constellations.} Conversely, let $P$
be a non-crossing partition of $\{1,2,\ldots,m\}$. We wish to find the
\gf\ of the $m$-constellations associated with $P$.

Take first a root $m$-gon $R$, and label it canonically, the
root-vertex being labeled $m$. Then, for each block $B$ of $P$, take a
constellation $C_B$ of outer-degree $md$, for some $d\ge 0$.
%with $d \ge |B|-1$. 
Label its root-vertex by $\min B$, and the other vertices in a canonical way. 

For each block $B$, we need to glue the component $C_B$ to the $m$-gon
$R$, and to keep track of the number of edges of $C_B$ that will
contribute to the outer-degree of the final constellation $C$. Let
$B=\{i_1, i_2, \ldots, i_k\}$ with $1\le i_1< i_2< \cdots< i_k\le
m$. If one walks around the root-face of $C_B$, starting from the
root-edge, the labels read at the corners of the root-face form the
word $u= (i_1 \cdots 
i_2 \cdots i_k \cdots)^d$. From now on, we identify the corners of
this face with the letters of $u$.
 For $r=1, \ldots , k$,  glue the $j_r$th corner labeled $i_r$ to
 the (unique) vertex labeled $i_r$ in $R$. To be consistent with the way
 we have chosen to root the components in the decomposition of a
 constellation, $j_1$ must be $1$. The condition for the final map to be
 planar is
$$
1\le j_k < \cdots < j_2 \le d
$$
(see Figure~\ref{fig:decomp-const}). Hence the component $C_B$ must
have outer-degree at least $m(k-1)$ and there are ${d\choose {k-1}}$ ways of
gluing $C_B$ to $R$.  

If $B$ is an internal block, none of its edges contribute to the
outer-degree of $C$. Otherwise,
$$
\delta(C_B,C) = \left\{
\begin{array}{ll}
md & \hbox{if } k=1,\\
i_k-i_1+m(j_k-1)& \hbox{if } k\ge 2.
\end{array}
\right.
$$
For $j_k$ fixed, the number of ways of choosing $ j_{k-1}, \ldots ,
j_2 $ is ${{d-j_k}\choose {k-2}}$.

By~\Ref{degree}, the outer-degree of the final constellation $C$ is thus
$$
\hbox{ext}(P)+ \sum_{B \hbox{ external}} \big( \max B -\min B + m \left(
j(B,C)-1\right)\big)
= m + m \sum_{B \hbox{ external}}\big( j(B,C)-1\big),
$$
where $j(B,C)=d+1$ if $B$ is a singleton and $C_B$ has outer-degree
$md$, and $j(B,C)$ is the number $j_k$ defined above if $B$ has at
least two elements.

Putting together the above results, one can write the \gf\ of
$m$-constellations associated with the partition $P$ as
$$
tu 
\prod_{k= 1}^{m-1} \left( \sum_{d\ge k-1}{d\choose k-1}F_d\right)^{I_k(P)}
\left( \sum_{d\ge 0} F_d u^d\right)^{E_1(P)}
\prod_{k= 2}^m \left( \sum_{d\ge k-1}\left(F_d \sum_{j=1}^d u^{j-1}
     {{d-j}\choose k-2}\right)\right)^{E_k(P)} ,
$$
where $F_d$ is the coefficient of $u^d$ in $F(u)$, that is, the \gf\
of constellations having outer-degree $md$, and $I_k(P)$
(resp. $E_k(P)$) denotes the number of internal 
(resp. external) blocks of cardinality $k$ in the partition $P$.
%the binomail coefficient ${{d-j}\choose k-2}$ corresponds to the
%choice of $j_{k-1}, \ldots , j_2$. 

Clearly, with the notation defined in the proposition,
$$
\sum_{d\ge k-1}{d\choose k-1}F_d = G_{k-1} \quad \hbox{and}
\quad \sum_{d\ge 0} F_d u^d=F(u).
$$
Now
\begin{eqnarray*}
\sum_{j=1}^d u^{j-1} {{d-j}\choose k-2} &=&
\sum_{i=0}^{d-1} u^{d-i-1} {i\choose k-2}\\
&=& \frac{u^{d-k+1}}{(k-2)!} \left. \frac{d^{k-2}}{dv^{k-2}}
\left(\frac{1-v^d}{1-v}\right)
\right|_{v=1/u} \hskip 10mm \hbox{(use Leibnitz' formula)}\\
&=& \frac{u^{d-k+1}}{(k-2)!} \left. \left( \frac{1-v^d}{(1-v)^{k-1}}(k-2)!
- \sum_{i=1}^{k-2} {{k-2}\choose i} {d\choose i}
\frac{i! (k-2-i)!v^{d-i}}{(1-v)^{k-i-1}}\right)\right|_{v=1/u}\\
&=& \frac{u^{d-k+1}}{(k-2)!}\left( (k-2)! \frac{u^d-1}{(u-1)^{k-1}}
u^{k-1-d}
-\sum_{i=1}^{k-2} (k-2)! {d\choose i}
\frac{u^{k-d-1}}{(u-1)^{k-i-1}}\right)\\
&=& \frac 1 {(u-1)^{k-1}} \left( u^d-1 - \sum_{i=1}^{k-2}{d\choose
  i}(u-1)^i\right)\\
&=&\frac 1 {(u-1)^{k-1}} \left( u^d - \sum_{i=0}^{k-2}{d\choose
  i}(u-1)^i\right).
\end{eqnarray*}
Consequently,
$$
\sum_{d\ge k-1} \left(F_d \sum_{j=1}^d u^{j-1} {{d-j}\choose k-2}\right)
=
\frac 1{(u-1)^{k-1}} \left( F(u) - \sum_{i=0}^{k-2} \frac{(u-1)^i}{i!}
      {F^{(i)}(1)}\right),
$$
and the proposition follows.
\cqfd

\noindent {\bf Note.} The above functional equations for constellations
were obtained a few years ago by the first author of this
paper. They were  used to conjecture that the
number of $m$-constellations having $n$ black faces is
$$
C_m(n)=
\frac{(m+1)m^{n-1}}{[(m-1)n+2][(m-1)n+1]}{mn\choose n}.$$
This conjecture was then proved  in a
bijective way~\cite{mbm-schaeffer-constellations}.

%%%%%%%%%%%%%%%%%%%%%%%%%%%%%%%%%%%%%%%%%%%%%%%%%
\subsection{Hard particles on planar maps}
\label{section-hard-equation}
%%%%%%%%%%%%%%%%%%%%%%%%%%%%%%%%%%%%%%%%%%%%%%%%%
We consider here rooted planar maps in which the vertices are either
vacant, or occupied by a particle, with the constraint that two
adjacent vertices cannot be both
occupied. In~\cite{mbm-schaeffer-ising}, it was shown 
that the \gf\ of such decorated maps (rooted at an edge with vacant
endpoints) is a specialization of the vertex-distribution \gf\ of
bipartite planar 
maps, and it was proved to be algebraic  as soon
as the degree of the vertices is bounded.

\begin{figure}[hbt]
\begin{center}
\input{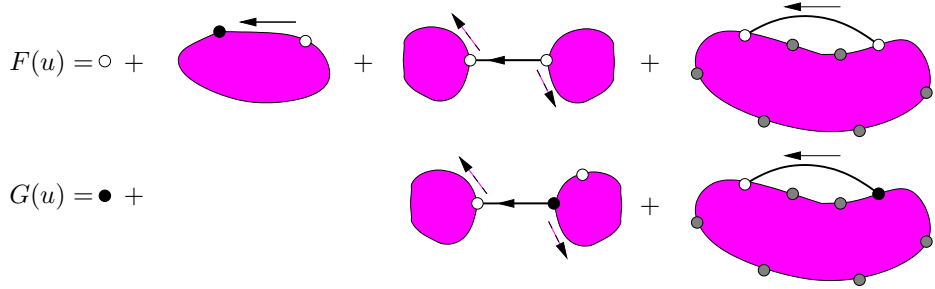}
\end{center}
\caption{The decomposition of planar maps carrying hard particles.}
\label{fig-hard}
%\hrule
\end{figure}

Here, we provide an independent approach for the case of unbounded
degrees. We say that an edge is \emm frustrated, if it has an occupied
endpoint (so that the other endpoint is vacant)\footnote{The
  terminology is standard in magnetism models like the Ising model.}. Let
$F(t,s,x,y,u)\equiv F(u)$ be
the \gf\ of maps with hard particles rooted at a \emm vacant, vertex,
counted by the number of edges ($t$),   frustrated
edges ($s$),   vacant vertices ($x$), occupied vertices ($y$), and
number of white corners in the infinite face ($u$).
Let $G(t,s,x,y,u)\equiv G(u)$ be defined similarly for maps with hard particles
rooted at an \emm occupied, vertex.  As observed by Gilles
Schaeffer~\cite{gilles}, it is not hard to adapt the equation written
for bipartite maps~\Ref{eq-bipartite} so as to obtain equations for
$F(u)$ and $G(u)$.  
\begin{Lemma}
\label{lemma-hard-particles}
The series $F(u)$ and $G(u)$ defined above are related by
\begin{eqnarray*}
F(u) &=& x -y + G(u) + tu^2 F(u)^2 + tu\, \frac{uF(u)-F(1)}{u-1},\\
G(u)&=& y+ tsu F(u) G(u) + tsu\, \frac{G(u)-G(1)}{u-1}.
\end{eqnarray*}
\end{Lemma}
\noindent
{\bf Proof.} As in Section~\ref{subsection-distribution}, these
equations follow from the deletion of the root-edge. From
Figure~\ref{fig-hard} one derives
$$
\begin{array}{lllllllllll}
F(u)&=& x &+& G(u)-y &+& tu^2F(u)^2 &+& t\sum_{j\ge 0} F_j \left( u+
\cdots + u^{j+1}\right), \\
G(u)&=& y &+& &+& tsuF(u)G(u) &+& ts\sum_{j\ge 0} G_j \left( u+
\cdots + u^{j}\right) ,
\end{array}
$$
where $F_j$ (resp. $G_j$) is the coefficient of $u^{j}$ in $F(u)$
(resp. $G(u)$). The result follows.
 \cqfd

Since the second equation is linear in $G(u)$, it is easy  to
eliminate $G(u)$. This gives a polynomial equation 
involving $F(u), F(1)$ and $G(1)$, and we can foresee that its
solution will be algebraic. We solve this equation in
Section~\ref{section-hard} (in the case $x=y=1$).

%%%%%%%%%%%%%%%%%%%%%%%%%%%%%%%%%%%%%%%%%%%%%%%%%
\section{From $3k$ to $2k$ equations: the role of the discriminant}
\label{section-discriminant}
%%%%%%%%%%%%%%%%%%%%%%%%%%%%%%%%%%%%%%%%%%%%%%%%%
We assume  again that $k+1$ power series in $t$, denoted $F(u), F_1,
\ldots, F_k$, are related by a functional equation of the form
\beq
P(F(u),F_1,\dots,F_k,t,u)=0.
\label{eq-main3}
\eeq
Here, $P(x_0, x_1, \ldots , x_k,t, v)$ is a polynomial with coefficients
in a field $\GK$, the $F_i$ belong to $\GK[[t]]$
and $F(u)$ belongs to $\GK[u][[t]]$. 
As discussed in Section~\ref{section-key}, for every fractional power
series $U\equiv U(t)$ such that
\beq
P'_0(F(U),F_1,\dots,F_k,t,U)=0,
\label{eq-derivee-x03}
\eeq
a system of three polynomial equations relating $U, F(U)$ and the unknown
functions $F_i$ holds:
$$
\left\{ 
\begin{array}{lll}
P\Big(F(U),F_1, \ldots , F_k, t,U\Big) &=&0,\\
P'_0 \Big(F(U),F_1, \ldots , F_k, t,U\Big)
&=&0,\\
P'_v\Big(F(U),F_1, \ldots , F_k, t,U\Big)
&=&0.
\end{array}
\right.$$
We say that the functional equation is \emm generic, if there exist
$k$ distinct series $U_i$ in $\overline \GK ^\f [[t]]$
satisfying~\Ref{eq-derivee-x03}. In this case, the strategy of
Section~\ref{section-key} provides a system of $3k$ polynomial 
equations relating the series $U_i, F(U_i)$ and $F_i$ for $1\le i \le
k$ (more precisely, a system of $3k+1$ equations, since one has to take into
account the fact that the $U_i$ are distinct, thanks to an equation
of the form~\Ref{Ui-distinct}).

The aim of this section is to eliminate the series $F(U_i)$, and to
reduce the system to $2k\,(+1)$ 
equations involving only the series $U_i$ and $F_i$. The key of this
reduction is the following theorem, which also considers the case
of multiple  roots $U_i$.

\begin{Theorem}\label{discriminant-original}
Assume that the functional equation~{\rm\Ref{eq-main3}} holds, and that the
series $U\in \overline \GK^\f[[t]]$ is a root of multiplicity $\ell$ of $P'_0 \left(F(u),F_1,
\ldots , F_k, t,u\right)$, with $\ell \ge 1$. 
Assume also that the degree of 
$P(x_0, \ldots , x_k, t,v)$ in $x_0$ is at least $2$, and  let
$\Delta(x_1, \ldots , x_k,t,v)$  be the discriminant of $P(x_0, \ldots
, x_k,t,v)$ with 
respect to $x_0$. Then, as a polynomial in $v$, $\Delta(F_1, \ldots ,
F_k,t,v)$ admits the series $U$ as a  root of multiplicity at
least $2\ell$. In other words, for $0\le i \le 2\ell-1$,
$$
\frac{\partial ^{i} \Delta}{\partial v^i}(F_1, \ldots , F_k,t, U)=0.
$$
\end{Theorem}
Recall that the discriminant of a polynomial $P(x)=a_nx^n+ \cdots +
 a_0$ such that $a_n\not = 0$ can be expressed as 
\beq
\label{discrim-def}
\Delta=(-1)^{n(n-1)/2} 
\left|
\begin{array}{ccccccccccccc}
1      &a_{n-1} & \cdots         &  a_2    &   a_1    & a_0    &   & \\
         &a_n     &         \cdots &         &   a_2      & a_1     & a_0 &\\
         &        &\ddots         &         &        & \ddots &     & \ddots\\
         &        &                &  a_n    &  &   \cdots     &     && a_0\\
n    & (n-1)a_{n-1} &\cdots     &  2a_2   &   a_1  & 0  & \\
         & na_n   &  \cdots     &       &      &a_1  &     & \\
&&\ddots&&&& \ddots\\
%&&\hskip 15mm \ddots&&&&&\ddots\\
         &        && & na_n  &        &  \cdots &        &        a_1
\end{array}
\right|.
\eeq
The above square matrix has size $2n-1$, and the coefficients that are not
indicated equal $0$.
In the generic case,  Theorem~\ref{discriminant-original} provides a
system of $2k$ equations:  
\beq
\label{system-2k}
\forall i \in [1,k], \hskip 8mm \left\{ 
\begin{array}{lll}
\Delta\left(F_1, \ldots , F_k, t,U_i\right) &=&0,\\
\Delta'_v\left(F_1, \ldots , F_k, t,U_i\right)
&=&0,
\end{array}
\right.
\eeq
which we complete with the distinctness
condition~\Ref{Ui-distinct}. But it may also happen that the series
$P'_0 \left(F(u),F_1,
\ldots , F_k, t,u\right)$ has a multiple root. In this case, the
system  derived from Theorem~\ref{discriminant-original} contains
more equations than unknowns. An example is 
provided in Section~\ref{section-multiple}.

In order to simplify  the proof of
 Theorem~\ref{discriminant-original}, 
 we first reduce  it to the case $U=0$.  Define 
$$
S(x,v):= P(x, F_1, \ldots , F_k, U+v) \quad  \hbox{and}\quad 
G(u):=F(u+U).$$
Then $S(x,v)$ is a polynomial in $x$ and $v$ with coefficients in
$ \GL=\overline\GK^\f((t))$, and $G(u)$ is a series of $\overline
\GK[u]^\f [[t]]$, and hence of $\GL[[u]]$. The functional
equation~\Ref{eq-main3} and the
assumption of  Theorem~\ref{discriminant-original} respectively imply
\beq
\label{assumptions-translated}
S(G(u),u) =0 \quad \quad \hbox{and} \quad\quad 
\frac{\partial S}{\partial x}(G(u),u) =u^\ell\Phi(u)
\eeq
with $\Phi(u) \in \overline \GK[u]^\f [[t]]\subset \GL[[u]]$ (the
second identity follows from Lemma~\ref{lemma-factorization}). Thanks
to this reduction, we will derive Theorem~\ref{discriminant-original}
from the following proposition. 

\begin{Propo}
\label{discriminant}
Let $\GL$ be an algebraically closed field,
% of characteristic $0$, 
and let $S(x,v)$ be a
polynomial  in $x$ with coefficients in $\GL^\f[[v]]$,
of degree $n\geq2$ in $x$.  Suppose that there exist two elements 
$G(u)$ and $\Phi(u)$ in $ \GL^\f[[u]]$ such that
$$
S(G(u),u) =0 \quad \quad \hbox{and} \quad\quad 
\frac{\partial S}{\partial x}(G(u),u) =u^\ell\Phi(u).
$$
Then the discriminant of $S(x,v)$ with respect to $x$, denoted
$\Delta(v)$, is
divisible by $v^{2\ell}$ in $\GL^\f[[v]]$.
\end{Propo}
The first step in the proof of Proposition~\ref{discriminant} is the
following ``exchange'' lemma. 
\begin{Lemma}
\label{lemma-discriminant}
Under the assumptions of Proposition~{\em\ref{discriminant}}, suppose,
moreover, that  $\frac{\partial^2S}{\partial x^2}(G(0),0)\ne0$. 
%  Let $\overline \GL$ be an algebraic closure of $\GL$. 
Then there exists $H(u)$ and $\Psi(u)$ in $  \GL ^\f[[u]]$ such that
$$
S(H(u),u) =u^{2\ell }\Psi(u) \quad \quad \hbox{and}\quad \quad 
\frac{\partial S}{\partial x}(H(u),u) =0.
$$
\end{Lemma}
\noindent{\bf Proof of Lemma~\ref{lemma-discriminant}.} We look for a
solution of the equation $S'_x(H(u),u)=0$
in the form $H(u)=G(u)+u^\ell Y(u)$, with $Y(u)\in  \GL^\f[[u]]$.
Using Taylor's formula, we write
\beq
\label{Taylor}
S(G(u)+z,u)=\sum_{i=1}^n\frac{z^i}{i!}\frac{\partial^iS}{\partial x^i}(G(u),u),
\eeq
and
$$
\frac{\partial S}{\partial x}(G(u)+z,u)=\sum_{j=0}^{n-1}\frac{z^j}{j!}
\frac{\partial^{j+1}S}{\partial x^{j+1}}(G(u),u).
$$
We thus want to find out whether there exists $Y\equiv Y(u)$ in
$ \GL^\f[[u]]$ satisfying
$$
\sum_{j=0}^{n-1}\frac{u^{\ell j}Y^j}{j!}
\frac{\partial^{j+1}S}{\partial x^{j+1}}(G(u),u)=0,
$$ 
that is,
$$
\frac{\partial S}{\partial x}(G(u),u)+\sum_{j=1}^{n-1}\frac{u^{\ell
    j}Y^j}{j!}\frac{\partial^{j+1}S}{\partial x^{j+1}}(G(u),u)=0,
$$
which, after dividing by $u^\ell$, reduces to
$$
\Phi(u)+\sum_{j=1}^{n-1}\frac{u^{\ell
    (j-1)}Y^j}{j!}\frac{\partial^{j+1}S}{\partial x^{j+1}}(G(u),u)=0.
$$
This is a polynomial equation in $Y$ with coefficients in $\GL^\f[[u]]$.
By  Theorem~\ref{thm-roots}, the number of roots lying in $ \GL^\f[[u]]$ is
$$
\text{deg}_Y\left(\Phi(0)+Y\frac{\partial^2S}{\partial
  x^2}(G(0),0)\right).
$$
The assumption of Lemma~\ref{lemma-discriminant} implies that this
degree is $1$, so 
that the equation $S'_x(H(u),u)=0$ admits a solution of the form
$H(u)=G(u)+u^\ell Y$, with $Y\equiv Y(u)$ in $ \GL^\f[[u]]$.
Then, by~\eqref{Taylor},
% if $H(u)=G(u)+u^\ell Y(u)$,
$$
\begin{array}{lll}
S(H(u),u) &=&\displaystyle \sum_{i=1}^n\frac{u^{\ell
 i}Y^i}{i!}\frac{\partial^iS}{\partial x^i}(G(u),u)\\ 
 &=&\displaystyle u^{2\ell }Y\Phi(u)+\sum_{i=2}^n\frac{u^{\ell
 i}Y^i}{i!}\frac{\partial^iS}{\partial x^i}(G(u),u), 
\end{array}
$$
which is divisible by $u^{2\ell }$ in $\GL^\f[[u]]$. 
\cqfd

\medskip

\noindent{\bf Proof of Proposition \ref{discriminant}.} 
Let 
$\tilde S(x,v)=S(x,v)+\eps(x-G(v))^2+\eps(x-G(v))^n$, 
where $\eps$ is a new indeterminate.
Then $\tilde S(x,v)$ belongs to $\GM^\f[[v]][x]$, where $\GM$ is the
algebraic closure of $\GL(\eps)$. Moreover,
$$
\tilde S(G(u),u)=0 \quad  \quad \hbox{and} \quad  \quad
\frac{\partial\tilde S}{\partial x}(G(u),u)=\frac{\partial S}{\partial
  x}(G(u),u)=u^\ell \Phi(u).
$$
Also,
$$
\frac{\partial^2\tilde S}{\partial
  x^2}(G(0),0)=\frac{\partial^2S}{\partial x^2}(G(0),0)+2\eps(1+
\delta_{n,2}) \ne0 \quad  \hbox{and} \quad  
\frac{\partial^n\tilde S}{\partial x^n}(x,0)=\frac{\partial^n\tilde
  S}{\partial x^n}(x,0)+ \eps ( 2\delta_{n,2}+ n!) \ne0.
$$
%For $\eps$ small enough, $\frac{\partial^2\tilde S}{\partial
%x^2}(G(0),0)\ne0$. 
 The discriminant of $\tilde S(x,v)$ with respect to $x$ can be written as
 follows~\cite[Ch.~V, \S~10]{lang}: 
\beq
\label{discrim-formula}
\tilde  \Delta(v)=\pm n^n a_n(v)^{n-1}\prod_{X(v)\in\R} \tilde S(X(v),v)
\eeq
where $a_n(v)\in \GM^\f[[v]]$ is the coefficient of $x^n$ in
$\tilde S(x,v)$ and $\R=\{X(v)\in\GM^\f((v))\ :\
\tilde S'_x(X(v),v)=0\}$.

 The condition $\frac{\partial^n\tilde S}{\partial x^n}(x,0)\ne0$,  combined
 with Theorem~\ref{thm-roots}, implies that all the elements of $\R$
are actually in $ \GM^\f[[v]]$.  Hence all the series $\tilde S(X(v),v)$, for
  $v\in \R$,  lie in  $\GM^\f[[v]]$. The condition
  $\frac{\partial^2\tilde S}{\partial 
  x^2}(G(0),0)\ne0$, combined with  Lemma~\ref{lemma-discriminant},
implies that one of the elements of $\R$, say $H(v)$, is such that
 $\tilde S(H(v),v)$ is 
divisible by $v^{2\ell}$. By~\eqref{discrim-formula},
$v^{2\ell }$ divides $\tilde \Delta(v)$ in $\GM^\f[[v]]$.

Since $\tilde \Delta(v)$ is a
polynomial in $\eps$,  this implies that each of its coefficients is
divisible by $v^{2\ell }$. Since its constant 
coefficient is equal to $\Delta(v)$, we conclude that $v^{2\ell }$ divides
$\Delta(v)$ in $\GL^\f[[v]]$.

\cqfd

\noindent{\bf Proof of Theorem \ref{discriminant-original}.} 
Let us return to~\eqref{assumptions-translated}. By
Proposition~\ref{discriminant}, \emm if the degree of, $S(x,v)$ \emm in, $x$
\emm is at least 2,,   the
discriminant of $S(x,v)$ with respect to $x$, denoted here
$\delta(v)$, has a root of multiplicity at least $2\ell$ at $v=0$.
We want to prove that the same holds for
 $\Delta(F_1, \ldots, F_k,t,U+v)$. How is this polynomial (in $v$) related
to $\delta(v)$?
\begin{itemize}
\item If $S(x,v)$ has degree $n$ in $x$, then $\Delta(F_1, \ldots,
F_k,t,U+v)=\delta(v)$, and the theorem follows from
Proposition~\ref{discriminant}. 
% it suffices to prove that $0$ is a root of
%$\delta(v)$ of multiplicity at least $2\ell$. 
% But $S(x,v)$ may have degree less that $n$ in $x$. 
\item If $S(x,v)$ 
has degree at most $n-2$, then~\eqref{discrim-def} shows that
$\Delta(F_1, \ldots, F_k,t,v)= 0$, and the result is trivial.
% re is nothing to prove. 
%If $n=2$, then $S(x,v)$ cannot have degree $1$ in $x$, 
\item 
If $S(x,v)$ has degree $n-1$, then
% subtracting from $P(x_0, \ldots, x_k,t,v)$ its leading monomial gives
%then there exists 
%a polynomial $Q(x_0, \ldots , x_k,t,v)$, of degree $n-1$ in $x_0$,
%such that $Q(F(u), F_1, \ldots, F_k,t,v)=0$ and $Q(x, F_1, \ldots,
%F_k,t,v)$ has degree $n-1$ in $x$. Moreover, 
\eqref{discrim-def} gives $\Delta(F_1, \ldots, F_k,t,U+v)= a_{n-1}(v)^2
\delta(v)$, where $a_{n-1}(v)$ is the coefficient of $x^{n-1}$ in
$P(x,F_1,\ldots, F_k,t,U+v)$. 
\begin{itemize}
\item If $n=2$, then $P(x, F_1, \ldots,
F_k,t,U+v)= a_1(v)x+a_0(v)$, where $a_0(v) $ and $a_1(v)$ belong to
$\overline\GK^\f[[t]][v]$. Then $P'_0(x,F_1, \ldots,
F_k,t,U+v)=a_1(v)$ and the assumption of
Theorem~\ref{discriminant-original} tells us that 
$0$ is a root of $a_1(v)$ of multiplicity $\ell$, and hence a root of
multiplicity $2\ell$ of $\Delta(F_1, \ldots, F_k,t,U+v)= a_{n-1}(v)^2$.
\item If $n\ge 3$, then $S(x,v)$ has degree at least $2$, and the
  theorem follows again from Proposition~\ref{discriminant}.  
\end{itemize}
\end{itemize}

% The equations~\eqref{assumptions-translated}, combined with
%  Proposition~\ref{discriminant}, imply that  $$\Delta(F_1, \ldots,
%  F_k,t,U+v)=v^{2\ell} \Psi(v)$$  where $\Psi(v) \in  \overline\GK^\f
%  ((t))[[v]]$. But $\Delta(F_1, \ldots, F_k,t,U+v)$ actually belongs
%  to $\overline \GK^\f[[t]][v]$, hence so does $\Psi(v)$. Replacing
%  $v$ by $v-U$ shows that $\Delta(F_1, \ldots, F_k,t,v)$ is divisible
%  by $(v-U)^{2\ell}$ in $\overline \GK^\f[[t]][v]$, which is the
%  expected result. 
\cqfd
%%%%%%%%%%%%%%%%%%%%%%%%%%%%%%%%%%%%%%%%%%%%%%%%%
\section{From $2k$ to $k$ equations: resultants and their generalization}
\label{section-resultants}
%%%%%%%%%%%%%%%%%%%%%%%%%%%%%%%%%%%%%%%%%%%%%%%%%
In the previous section, we have shown how to reduce our polynomial
system to $2k+1$ equations, at least in a generic situation.
%  case where  the $k$ series $U_i$ are distinct. 
%The first $2k$ equations are given
%by~\eqref{system-2k}, and the last one tells us that the $U_i$ are distinct.
%The whole set of equations can be described by saying that the
This system says that the polynomials (in $v$) $\Delta(F_1, \ldots ,
F_k,t,v)$ and 
$\Delta'_v(F_1, \ldots , F_k,t,v)$ have $k$ distinct roots in common.
It is well-known that two polynomials have one  root in common if and only
if their resultant vanishes. This is the result we generalize in this
section: we give a criterion
%two criteria 
that tells  when two polynomials $P$ and $Q$ have $k$ roots
in common. If the respective degrees of $P$ and $Q$ are $m$ and $n$,
%the first criterion involves many  determinants of order $m+n-2k+2$,
this criterion involves $k$ determinants 
of respective order $m+n$, $m+n-2$,\ldots, $m+n-2k+2$.
%, the determinant of order $m+n$ being the resultant of $P$ and $Q$. 
In a generic situation,
these determinants directly provide $k$ equations between the series
$F_1, \ldots, F_k$, with no mention of the series $U_i$. Whether these
equations are as small as they can be is another story...

\medskip
Let $P(X)=\sum_{i=0}^ma_iX^i$ and $Q(X)=\sum_{i=0}^nb_iX^i$, where the
coefficients $a_i$ and $b_i$ belong to a field $\GL$.
% , and $a_mb_n\not = 0$. 
For $0\le k<\min(m,n)$, we define a matrix $\mathcal
S_k(P,Q)$ having $m+n-2k$ rows and  columns by:
$$
\mathcal S_k(P,Q)=\left(
\begin{array}{ccccccccccc}
%\displaystyle
%
%
%   Premiere ligne
%
%
a_m      &        & \cdots &        &   a_0    & & & \\
%         &a_m     &        & \cdots &          &   a_0   &   &\\
         &\ddots  &        &          &         & \ddots & \\
         &        &  a_m   &          &\cdots     & &a_0   \\
&        &        &  \ddots   &          &&\vdots   &    \\
	 &         &        &       &   a_m     &     \cdots      & a_k\\
b_n      &        & \cdots &        &          &   b_0      &   & \\
         & \ddots &        &        &          &         &  & \\
         &        &  b_n   &        &  \cdots  &           &b_k    \\
%    &        &  b_n   &        &  \cdots  &   &        &b_k

%
\end{array}
\right)
$$
where the first $n-k$ rows are filled up with the coefficients
$a_i$ of $P$ and the $m-k$ last ones by the coefficients $b_i$
of $Q$. The other entries are zero. In particular, $\mathcal S_0(P,Q)$
is the \emm Sylvester matrix, of $P$ and $Q$ (and its determinant is the
resultant of $P$ and $Q$). In general, $\mathcal S_k(P,Q)$
is obtained by
% the matrix with $m+n-2k+2$ rows and $m+n-k+1$ columns obtained by
% cancelling 
deleting the $k$ last rows of $a$'s, the $k$ last rows of $b$'s,
%rows $n$,  $n-1$,..., $n-k+2$, $n+m$,..., $n+m-k+2$, 
and the $2k$ rightmost columns 
%$m+n$,..., $m+n-k+2$ 
in the Sylvester matrix of $P$ and $Q$. 

The following theorem is a simple adaptation
 of~\cite[Prop.~4.33]{basu}.

\begin{Theorem}\label{2k-kcritere2}
Let $k\le \min(m,n)$. If the polynomials $P$ and $Q$ have $k$ common roots,
counted with multiplicities, then  for $0\le i \le k-1$,
$$\det\,\mathcal S_{i}(P,Q)=0.$$
Conversely, if the above determinants vanish, then either  $P$ and $Q$
have $k$ common roots, or $a_m=b_n=0$.
\end{Theorem}

%%%%%%%%%%%%%%%%%%%%%%%%%%%%%%%%%%%%%%%%%%%%%%%%%
\section{A new proof of Brown's theorem}
\label{section-brown}
%%%%%%%%%%%%%%%%%%%%%%%%%%%%%%%%%%%%%%%%%%%%%%%%%
Theorem~\ref{thm-brown} below is essentially due to
Brown~\cite{brown-square}. It has been used several times in the past
to solve functional equations of the form~\Ref{eq-generic}. Its application is
straightforward for quadratic  equations
(see~\cite{bender-canfield,gao-5-connected} and the discussion at the
end of this section), but more elusive when the 
degree in $F(u)$ is larger (see~\cite[Section~4]{brown-square} for a
solution of the 
cubic equation of Section~\ref{section-dissections} based on this theorem).

In this section, we give a new proof of, and a new point of view on
 Brown's theorem\footnote{Or maybe we should write that \emm we give  a
  proof of, this theorem, since there seems to be a mistake in Brown's
  proof~\cite{brown-square}: in the equation that follows (2.12), why aren't
  there any terms $U_{r-1}V_1, \ldots , U_0V_r$?}. Here is our
formulation of this theorem.

\begin{Theorem}\label{thm-brown}
%Let $\GK$ be a field of characteristic $0$. 
%\mps{algebraically closed?}
Let $\Delta(t,u) \in\GK[[t]][u]$, where  $\GK$ is a field.
%  be a non-zero polynomial of $\GK[[t]][u]$. 
If $\Delta$ has a square root in $\GK[[t,u]]$, then it can be factored as
$$
\Delta(t,u)= c ^2 t^{2p} (1+tS(t)) \Big(1+tuR_1(t,u)\Big)
\Big(u^d+tR_2(t,u)\Big)^2 
\prod_{i=1}^k \left( \left(1-\frac
u{\alpha_i}\right)^{d_i}+tuQ_i(t,u)\right),
$$
where
\begin{itemize}
\item[--] $p,d$, and the $d_i$'s are nonnegative integers,
\item[--] $c $ belongs to $\GK$ and the $\alpha_i$'s belong to
  $\overline \GK$, the algebraic closure of $\GK$,
\item[--] $S(t)\in \GK[[t]]$, 
\item[--] $R_1(t,u)$ and $ R_2(t,u)$ belong to $ \GK[[t]][u]$,  with
$\deg_u(R_2) <d$,
\item[--]
  $Q_i(t,u)$ belongs to $ \ov \GK[[t]][u]$, and $\deg_u(Q_i)<d_i$. 
\end{itemize}
Moreover, if  $\Delta$ has a square root in
$\GK[u][[t]]$, then it can be factored as
$$
\Delta(t,u)= c^2 t^{2p} (1+tS(t)) \Big(1+tuR_1(t,u)\Big) \Big(u^d+tR_2(t,u)\Big)^2
\prod_{i=1}^k \left( \left(1-\frac
u{\alpha_i}\right)^{d_i}+tuQ_i(t,u)\right)^2,
$$
with the same conditions as above.
\end{Theorem}
What is remarkable in the above factorizations is the fact that some
factors are 
squared. We will derive this theorem from the combination of two results. The
first one is a factorization theorem which has an
independent interest and  will be used in 
Sections~\ref{section-practise} and~\ref{section-distribution}. 
\begin{Theorem} [{\bf Factorization Theorem}]
\label{thm-factorization}
%Let $\GK$ be a field of characteristic $0$, and 
Let $\Delta(t,u)$ be a non-zero polynomial of
$\GK[[t]][u]$, where $\GK$ is a field.
%$\Delta(t,u) \in \GK[[t]][u]$. 
Then $\Delta$ admits a unique factorization as
$$
\Delta(t,u)= c t^{p} (1+tS(t)) \Big(1+tuR_1(t,u)\Big) \Big(u^d+tR_2(t,u)\Big)
\prod_{i=1}^k \left( \left(1-\frac
u{\alpha_i}\right)^{d_i}+tuQ_i(t,u)\right),
$$
with the same conditions as in Theorem~{\rm\ref{thm-brown}.}
The roots of $\Delta(t,u)$ that are infinite (resp. zero, finite and
non-zero) at $t=0$ are the roots of the first (resp. second, third)
factor above. 
\end{Theorem}
\noindent
{\bf Proof.} Let us first recall that the units of
$\GK[[t]][u]$ and $\GK[[t]]$ coincide, and are the series
$c(1+tS(t))$, where $c\in\GK\setminus\{0\}$ and $S(t)  \in
\GK[[t]]$. 

Now consider an irreducible polynomial  of
$\GK[[t]][u]$, denoted  $P(t,u)$, of degree $d$ in $u$. By definition,
$P$ is not a unit. If $d=0$, then
$$
P(t,u)=tI(t),
$$ 
where $I(t)$ is  unit of $\GK[[t]]$.
 If $d>0$, then $P(0,u)\not =0$ (otherwise
$P$ would be divisible  by $t$). Moreover, $P(t,u)$
is also  irreducible in $\GK((t))[u]$. The roots $U_1, \ldots, U_d$ of
$P$ are of the form~\cite[Prop.~6.1.6]{stanley-vol2} 
$$
U_i= \sum_{n\ge n_0} a_n \left( \xi^i t^{1/d}\right)^{n}
$$
where  $n_0 \in \zs\cup\{+\infty\}$,  $a_{n_0}\not =0$, the
coefficients  $a_n$ lie in $\ov \GK$, and $\xi$ is a primitive $d$th root
of unity in $\ov \GK$.
%  If $U_i=0$, we consider that $n_0=+\infty$, and we otherwise assume that. 
We consider three cases, depending on whether
$n_0$ is negative, positive or zero.

\smallskip
\noindent
{\bf Case 1: $n_0<0$.}
Then all the roots of $P$  are infinite at $t=0$. By
Theorem~\ref{thm-roots}, $\deg_uP(0,u)=0$. Thus $P$ can be written
$
P(t,u)=P(t,0)+tuR(t,u)
$
with $\deg_uR(t,u)=d-1$. Since by assumption $P(0,u)\not = 0$, we have
$P(0,0) \not =0$, so that $P(t,0)$ is a unit in $\GK[[t]]$. Denoting 
$P(t,0)=P_0(t)$, we have
$$
P(t,u)=P_0(t)\left( 1 +tu R_1(t,u)\right)
$$ 
with $R_1(t,u)\in \GK[[t]][u]$.

\smallskip
\noindent
{\bf Case 2: $n_0>0$.} All the roots  of $P$ are zero at $t=0$.  By
Theorem~\ref{thm-roots}, $\deg_uP(0,u) =d$.
 More precisely, $P(0,u)=c u^d$ for some $c\in \GK\setminus\{0\}$.
Denoting by $P_d(t)$ the
coefficient of $u^d$ in $P$, we thus have  $P(t,u)=P_d(t) u^d +
tR(t,u)$ where $\deg_uR(t,u) <d$, and $P_d(0)=c \not =0$. Thus $P_d(t)$
is  a unit of $\GK[[t]]$, and we can write
 $$
P(t,u)=P_d(t) \Big( u^d+tR_2(t,u)\Big),
$$ 
where  $R_2(t,u)\in \GK[[t]][u]$ and $\deg_uR_2(t,u) <d$.

\smallskip
\noindent
{\bf Case 3: $n_0=0$.} All the roots  of $P$ are equal to some $\alpha \not =
0$ when $t=0$, with $\alpha\in\overline \GK$. As in Case 2,  $\deg_uP(0,u)
=d$, and more precisely 
$$
P(0,u)=c\left(1-\frac u \alpha\right)^d 
$$
where $c\in \GK\setminus\{0\}$. In particular, $P(0,0)=c\not =0$, so
that $P(t,0)\equiv P_0(t)$ is a unit in $\GK[[t]]$. Thus we can write
$P(t,u)=P_0(t)\left( 1+ uR(t,u)\right)
$
where $R(t,u) \in \GK[[t]][u]$ and $\deg_uR(t,u)=d-1$. Setting $t=0$ gives
$
P(0,u) = P_0(0)\left( 1+ uR(0,u)\right)= P_0(0)\left(1-\frac u
\alpha\right)^d 
$.
Finally,
$$
P(t,u)=P_0(t)\left( \left(1-\frac u \alpha\right)^d +tu
Q(t,u)\right)
$$
where
$$
Q(t,u)= \frac{R(t,u)-R(0,u)} t
$$
belongs to $\GK[[t]][u]$ and has degree at most $d-1$ in $u$.

\medskip
Now take $\Delta \in  \GK[[t]][u]$, as in the statement of the
theorem. Factor $\Delta$ into irreducible polynomials of $
\GK[[t]][u]$. Write each irreducible factor in the above canonical
form. Then, group together the irreducible factors whose roots are
infinite (resp. zero, equal to $\alpha_i\not = 0$) at $t=0$. This gives for
$\Delta$ a factorization of the prescribed form.
%$$
%\Delta(t,u)= ct^p\left(1+tS(t)\right) \Big(1+tuR_1(t,u)\Big) \Big(u^d+tR_2(t,u)\Big)
%\prod_{i=1}^k \left( \left(1-\frac
%u{\alpha_i}\right)^{d_i}+tuQ_i(t,u)\right),
%$$
%where 
%$C(t) \in \GK((t))$ and 
%the polynomials $R_1, R_2 $ and $Q_i$
%satisfy the required conditions. 
% Finally, note that the valuation of $\Delta$ in $t$ equals the
% valuation of $C(t)$, which is, subsequently, nonnegative. 
% Writing $C(t)=t^p(1+tS(t))$ 

The uniqueness of this factorization is a consequence of  the two
 following facts:

-- the roots of the first (resp.~second, third) factor
%  of positive degree in $u$ 
are exactly the roots of $\Delta(t,u)$ that are infinite
(resp.~zero, equal to $\alpha_i \not = 0$) at $t=0$, 

--  these
factors are normalized (they have either constant term 1, or leading
coefficient $1$). 

\noindent This concludes the proof.
\cqfd

In order to prove Brown's theorem (Theorem~\ref{thm-brown}), we only
need to combine the above factorization theorem with the following
proposition.
\begin{Propo}\label{propo-even-mult}
Let $\Delta(t,u) \in \GK^\f[[t]][u]$, where $\GK$ is a field.
% of characteristic $0$.  
The roots of $\Delta(t,\cdot)$ belong to $\overline
\GK ^\f ((t))$. 

If $\Delta$ has a square root in
$\GK[[u]]^\f[[t]]$, then every root $U$ of $\Delta$ that vanishes at $t=0$
has an even multiplicity in $\Delta$. 

If $\Delta$ has a square root in
$\GK[u]^\f[[t]]$, then every root $U$ of $\Delta$ that is finite at $t=0$
has an even multiplicity. 
\end{Propo}
\noindent
{\bf Proof.} 
% The proofs of both statements are very similar, and we only prove
% the first one. 
%
Assume $\Delta(t,u) = \delta(t,u)^2$, with $\delta \in
\GK[[u]]^\f[[t]]$. Let $U\equiv U(t)$ be a root of $\Delta$ that
vanishes at $t=0$. Then $\delta(t,U)$ is a well-defined series in
$\GK^\f[[t]]$, which must be $0$. Thus $U$
 is a root of $\delta(t,u)$, and by Lemma~\ref{lemma-factorization},
\beq
\label{Dd}
\delta(t,u)= (u-U) \Psi(t,u)
\eeq
where $\Psi(t,u) \in \overline\GK[[u]]^\f[[t]]$, so that
$$
\Delta(t,u)= (u-U)^2 \Psi(t,u)^2.
$$
Thus $U$ is a root of $\Delta$ of multiplicity at least $2$.

More generally, let us prove by induction on $m\ge0$ that, if
$U$ has multiplicity at least $2m+1$ in 
$\Delta$, then it actually has multiplicity at least $2m+2$. The case
$m=0$ has just been proved. Now take $m\ge1$, and assume
$$
\Delta(t,u)=(u-U)^{2m+1}\Delta_1(t,u)=\delta(t,u)^2
$$
with $\Delta_1(t,u) \in \overline\GK^\f[[t]][u]$. As argued above, $U$ is a
root of $\delta(t,u)$, and the factorization~\Ref{Dd} gives
$$
\tilde\Delta(t,u):=(u-U)^{2m-1}\Delta_1(t,u)=\Psi(t,u)^2.
$$
 The induction hypothesis implies that
$U$ is a root of $\tilde\Delta(t,u)$ of multiplicity  at least $2m$, and
 thus a  root of $\Delta(t,u)$ of multiplicity  at least $2m+2$.  This
 completes the proof of the first statement of the proposition. 

The proof of second statement  is very
similar. It relies on the fact that if $\delta(t,u)$  lies in
$\GK[u]^\f[[t]]$, then all 
roots $U$ of $\Delta$ that are finite at $t=0$ can be substituted for
$u$ in $\delta(t,u)$.
\cqfd

\noindent
We are finally ready for a \\
\noindent
{\bf Proof of Theorem~\ref{thm-brown}.}
Take $\Delta(t,u) \in \GK[[t]][u]$ and consider its canonical
factorization, given by Theorem~\ref{thm-factorization}.  Assume
$\Delta(t,u) = \delta(t,u)^2$, 
with $\delta \in \GK[[u,t]]$. If $q$ is the valuation in $t$ of
$\delta$, then  the valuation in $t$ of $\Delta$ is $p=2q$. Thus
 $p$ is even. Now 
$$
t^{-2q} \Delta(t,u) = \left(t^{-q}\delta(t,u)\right)^2.
$$
Setting $t=0$ in this identity shows that the constant $c$ occurring
in the canonical factorization of $\Delta$  is a square of $\GK$.

 By Proposition~\ref{propo-even-mult}, each root of $\Delta$ that
 vanishes at $t=0$ has an even multiplicity in $\Delta$. This means
 that every irreducible factor of $\Delta$
% (over $\GK((t))$) 
occurring
 in the term $(u^d+tR_2(t,u))$ actually occurs an even number of
 times. This implies that $d$ is even, and that this term can be
 factored  as $(u^{d/2}+t\tilde R_2(t,u))^2$.
  This completes the proof of the first statement.

The proof of the second statement is very similar: now, each root of
$\Delta$ that is finite at $t=0$ must have an even multiplicity.
\cqfd

\medskip

Let us finally discuss how Brown's theorem may be used to solve a
\emm quadratic, equation with one catalytic
 variable~\cite{bender-canfield,gao-5-connected}. We start from a
 $(k+1)$-tuple of series, denoted $F(u), F_1,
 \ldots , F_k$, such that $F(u)\in \GK[u][[t]]$ and $F_i \in\GK[[t]]$
 for all $i$. We assume they satisfy 
\beq
\label{eq-quadratic}
\Big(2 a F(u)+b\Big)^2  = \Delta(u),
\eeq
where $a,b$ and $\Delta$ are polynomials in $F_1, \ldots, F_k, t$
and $u$, with coefficients  $\GK$. Obviously, $\Delta$ has a
square root in $\GK[u][[t]]$ (namely, the series $2 a F(u)+b$). Hence
the second part of Theorem~\ref{thm-brown}, applies:
the canonical factorisation of $\Delta$ contains several squared factors.

Let us now adopt the notation of Theorem~\ref{thm-brown}.  In order to
determine the degrees in $u$ of $R_1, R_2$ and the $Q_i$, one has to 
decide how many  roots of $\Delta$ are infinite (resp. equal to
zero, equal to $\alpha_i$) when $t=0$. This can be done routinely
using Newton's polygon method. (Curiously, these degrees are only \emm
guessed, in~\cite{bender-canfield} and~\cite{gao-5-connected}. This
 forces the authors to check afterwards the validity of their
 assumption.) One then introduces a new set of indeterminates (the
 coefficients of the polynomials $R_1, R_2$ and  $Q_i$) and obtains a
 system of polynomial equations by comparing the coefficient of $u^j$
 in $\Delta$ and in its factorisation, for all $j$. This is
 illustrated in 
Sections~\ref{section-factorization} and~\ref{section-distribution}
 (even though we do not use Brown's theorem, but rather
a combination of our general strategy with the factorization theorem,
Theorem~\ref{thm-factorization}). 

To conclude, let us underline one important difference between the
quadratic case and the general case. As shown by~\eqref{eq-quadratic},
in the quadratic case, \emm every, root of $\Delta$ that is finite at $t=0$
% belonging to $\GK^\f[[t]]$ 
has an even multiplicity. 
For a general (i.e., non-quadratic) equation,
Theorem~\ref{discriminant-original} exhibits  a 
certain number of multiple roots of $\Delta$, which are finite at
$t=0$.  But $\Delta$ may also have simple 
roots (or roots of odd multiplicity) that are
finite at $t=0$. In the example of
Section~\ref{section-hard} below, $\Delta$ has two simple roots that
are finite at $t=0$.
%, but only two of them have an even multiplicity.

%%%%%%%%%%%%%%%%%%%%%%%%%%%%%%%%%%%%%%%%%%%%%%%%%
\section{Practical strategies}
\label{section-practise}
%%%%%%%%%%%%%%%%%%%%%%%%%%%%%%%%%%%%%%%%%%%%%%%%%
The general strategy presented in Section~\ref{section-key} to solve functional
equations of the form $P(F(u), F_1, \ldots , F_k,t,u)=0$ yields a
system of polynomial equations ($3k+1$ equations in a generic
 case) relating the unknown series $F_i$, some
auxilliary series $U_i$, and the values of
$F(U_i)$. Section~\ref{section-discriminant} 
performs the elimination of the
$F(U_i)$, yielding a system of $2k+1$
equations. Section~\ref{section-resultants} even goes further by 
eliminating the $U_i$, 
but the $k$ equations it provides are often, in practise,  unnecessarily
big. At any rate, it is always
easy to write a system of  $2k+1$ equations relating the series $F_i$
and $U_i$. 

In most combinatorial problems, one is interested in finding the
minimal equation satisfied by $F_1$, or at least a ``nice'' system
involving all the $F_i$, if such a system exists. 
 As discussed
in~\cite[Section 4]{fla-poly-alg}, three main methods can be used to reduce
further the size of our  system:  the
paper-and-pencil approach, the resultant approach, and the Gröbner
basis approach. Note that our system of  $2k+1$ equations contains $k$
  times the ``same'' pair of equations (see~\eqref{system-2k}), which
  means that the   elimination of the $U_i$ must be performed with
  care not to   loose any information.

The paper and-pencil approach has been amply illustrated in
 Section~\ref{section-examples}. In almost all examples  presented
 there, there was actually a \emm single, unknown function $F_1$. In
 this  case, as soon as one finds a series $U$ that cancels $P'_0(F(u),
 \ldots, t,u)$, the discriminant $\Delta$ has a double root
 (Theorem~\ref{discriminant-original}), and one obtains immediately an
 equation for $F_1$ by writing  that 
\begin{center}
\emm the discriminant of the   discriminant vanishes,.
\end{center} 
We also studied in Section~\ref{section-examples} one equation involving two
unknown series $F_i$,
%  that was considered in% Section~\ref{section-examples} 
but it was linear   in $F(u)$ and of an especially simple form.

In this section, we gather a number of practical strategies that
 permit to solve bigger examples.   We advise the reader who would be
 interested in the practical aspects of our method to read what follows
 with a computer  algebra system at hand. All the strategies we suggest have
 been tested on the same example (except the Gröbner one, which seems
 to be  too  brutal to work). Two more 
 examples are provided in Sections~\ref{section-distribution}
 and~\ref{section-hard}. 
%%%%%%%%%%%%%%%%%%%%%%%%%%%%%%%%%%%%%%%%%%%%%%%%%%%%%%%%%%%%%%%%%%%%
\subsection{Brute force on $3k+1$ equations}
The laziest approach naturally consists in feeding a Gröbner basis package with
the $3k(+1)$ equations obtained in the generic case, and let it
work. The aim is to obtain either a polynomial system defining the
series $F_i$, or 
%to eliminate most of them to obtain 
a single algebraic equation for, say, $F_1$. One has to choose carefully a
monomial order. See~\cite{cox} for generalities on Gröbner bases,
and~\cite{salvy} for a recent study of the complexity of Gröbner
computations.

Unfortunately, this lazy approach often fails, because the computation
tends to take forever. This is why we only give here a very simple ---
and somewhat degenerate --- example. 

Return to the second example of Section~\ref{section-kernel},
Eq.~\eqref{+3-2}. Form a set $S$ of 5 equations consisting
of~\eqref{kernel2-3} for $i=1,2$, the right-hand side
of~\eqref{complet2-3} for $i=1,2$ 
again, and the distinctness condition 
$X(U_1-U_2)=1$. The Maple command\\
{\tt Groebner[univpoly]}$(F_1,S,\{X,U_1,U_2,F_1,F_2\})$ directly gives
$$
{F_1}=1+2\,{t}^{5}{{F_1}}^{5}-{t}^{5}{{F_1}}^{6}+{t}^{5}{{
F_1}}^{7}+{t}^{10}{{F_1}}^{10}.
$$

%%%%%%%%%%%%%%%%%%%%%%%%%%%%%%%%%%%%%%%%%%%%%%%%%
\subsection{Bare hands elimination on $2k+1$ equations}

\subsubsection{The number of $3$-constellations}
\label{section-3const}
Let us consider now the equation~\eqref{eq-3const} that defines the \gf\ of
3-constellations. It has degree 3 in $F(u)$ and contains two unknown
series $F_1=F(1)$ and $F_2=F'(1)$. Multiplying by $(u-1)^2$ gives
an equation of the form $P(F(u), F_1,
F_2,t,u)=0$. Theorem~\ref{thm-roots}, applied to 
$P'_0(F(u),F_1,F_2,t,u)$, shows that this series  has two
roots, $U_1$ and $U_2$. Indeed, $P'_0(F(u),F_1,F_2,t,u)$ reduces to
$(u-1)^2$ when $t=0$. The original functional equation gives the first
terms of $F(u)$:  
$$
F(u)= 1+tu+3 \left( u+1 \right) u{t}^{2}+2 \left( 6\,{u}^{2}+10\,u+11
 \right) u{t}^{3} +O(t^4),
$$
and the equation $P'_0(F(U_i),F_1,F_2,t,U_i)=0$ provides the first terms of
$U_1$ and $U_2$:
$$
U_{1,2}=1\pm t^{1/2} +2\,{t}\pm 5\,{t}^{3/2}+15\,{t}^{2}\pm
48\,{t}^{5/2} +O(t^3).
$$
In particular, the series $U_i$ are distinct.
Let $\Delta(F_1,F_2,t,v)\equiv \Delta(v)$ be the discriminant of
$P(x,F_1,F_2,t,v)$, 
taken with respect to $x$. By Theorem~\ref{discriminant-original},
$\Delta(v)$ admits $U_1$ and $U_2$ as multiple roots. But $\Delta(v)$
factors as $tv(v-1)^4R_1(v)$, where $R_1$ is a 
polynomial in $F_1,F_2, t$ and $v$  of degree 5 in $v$. Since $U_i\not
=0$ and $U_i\not =1$, we conclude 
that $U_1$ and $U_2$ are double roots of $R_1$.
Let $R_2(v)$ be the derivative of $R_1$ with respect to $v$. Then
$R_1$ and $R_2$  have the roots $U_1$ and $U_2$ in common.

\begin{figure}[hbt]
\begin{center}
\input{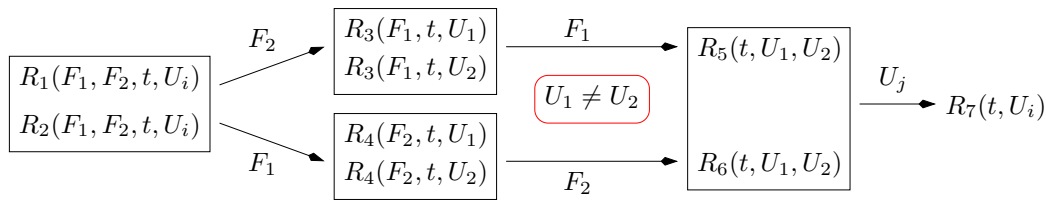}
\end{center}
\caption{The bare-hand elimination procedure for $k=2$.}
\label{fig-elimination}
%\hrule
\end{figure}

The rest of the elimination procedure is schematized in
Figure~\ref{fig-elimination}. The labels on the arrows indicate which
variable is eliminated (using a resultant) at each stage.
%For the sake of simplicity, we do not
%mention $t$ in the variables of the polynomials.

We first eliminate $F_2$ between $R_1(F_1,F_2,t,U_i)$ and
$R_2(F_1,F_2,t,U_i)$.  The polynomial thus obtained
factors in much smaller terms than $R_1$ or $R_2$. Knowing
the first terms of $F(u), U_1$ and $U_2$ allows us to decide which
factor vanishes.
% when $v=U_i$. 
We thus obtain $R_3(F_1,t,U_i)=0$, with
$$
R_3(F_1,t,v)= {t}^{2}{v}^{4}{F_1}^{2}
+2\,t{v}^{2}F_1(v-1)(v-3)
%-8\,t{v}^{3}F_1+6\,t{v}^{2}F_1
-4\,t{v}^{2}+{v}^{4}-8\,{v}^{3}+22\,{v}^{2}-24\,v+9.
$$
Similarly, if we eliminate $F_1$ between $R_1(F_1,F_2,t,U_i)$ and
$R_2(F_1,F_2,t,U_i)$, and then choose the right factor, we obtain
an equation of the form $R_4(F_2,t,U_i)=0$, of degree $8$ in $U_i$.
% with
%$$
%R_4(F_2,t,v)={t}^{4}{v}^{8}{F_2}^{2}-2\,{t}^{2}{v}^{4} \left( t{v}^{4}-{v}^{4
%}+8\,{v}^{3}-3\,t{v}^{2}-22\,{v}^{2}+24\,v-9 \right) F_2
%$$
%$$ +81+972
%\,{v}^{2}-6\,{v}^{6}{t}^{2}-1200\,{v}^{3}+886\,{v}^{4}-46\,t{v}^{2}-
%400\,{v}^{5}+108\,{v}^{6}+{v}^{8}-16\,{v}^{7}-162\,t{v}^{4}+{t}^{2}{v}
%^{8}-2\,{v}^{8}t+9\,{t}^{2}{v}^{4}-54\,t{v}^{6}+112\,{v}^{5}t-432\,v+
%136\,t{v}^{3}+16\,{v}^{7}t.
%$$

We now eliminate $F_1$ between $R_3(F_1,t,U_1)$ and
$R_3(F_1,t,U_2)$. The resultant naturally contains a factor
$(U_1-U_2)$, \emm which we know to be non-zero,. Choosing the right factor
among the remaining ones provides a first equation between $U_1$ and
$U_2$, of the form $R_5(t,U_1,U_2)=0$. Similarly, eliminating $F_2$
between $R_4(F_2,t,U_1)$ and 
$R_4(F_2,t,U_2)$ provides another such equation, say
$R_6(t,U_1,U_2)=0$. Finally, eliminating 
one of the $U_i$'s between $R_5$ and $R_6$ gives $R_7(t,U_i)=0$, with
$$
R_7(t,v)=
\left( t-4 \right) ^{3}{v}^{6}-4\, \left( 21\,t+44 \right) 
 \left( t-4 \right) {v}^{5}- \left( 180\,t+2944+27\,{t}^{2} \right) {v
}^{4}
$$
$$
-18\, \left( -332+15\,t \right) {v}^{3}+27\, \left( -235+9\,t
 \right) {v}^{2}+3402\,v-729.
$$
We have finally obtained the algebraic equation (on $\qs(t)$, of degree 6) satisfied
by each of the series $U_i$. It remains to eliminate $U_1$ between
$R_3(F_1,t,U_1)$ and $R_7(t,U_1)$ to obtain, by extraction of the
relevant factor,  the (cubic) algebraic equation satisfied by $F_1$:
\beq
\label{3const-sol}
F_1=1-47t+3t^2 +3t(22-9t)F_1+9t(9t-2)F_1^2-81t^2F_1^3.
\eeq
Recall that $F_1$ counts $3$-constellations by their number of black
triangles.

%%%%%%%%%%%%%%%%%%%%%%%%%%
\subsubsection{An example with multiple roots $U_i$}
\label{section-multiple}
%%%%%%%%%%%%%%%%%%%%%%%%%%%%%%%%%%%%%%
Consider  the functional equation
\begin{equation}\label{example-doubleroot}
F(u)=u+t\left(F(u)^3-3+2\,\frac{F(u)-F(0)}u- 
t\,\frac{F(u)-F(0)-uF'(0)}{u^2}\right).
\end{equation}
Clearly, it has a unique power series solution. 
The first terms of the expansion of $F$ are:
$$
F(u)=u+(u^3-1)t+u^2(3u^3-1)t^2+3u^4(4u^3-1)t^3+u^6(55u^3-12)t^4+\cdots
$$
After multiplying by $u^2$,  our functional equation reads
$P(F(u),F(0),F'(0),t,u)=0$, for some polynomial $
P(x_0,x_1,x_2,t,u)$.
%=-u^3+3tu^2-tx_0^3u^2+(u-t)^2x_0+2tux_1-t^2x_1-t^2ux_2.
%
We are looking for fractional series $U$ that satisfy
$P'_0(F(U),F(0),F'(0),t,U)=0$, that is
\begin{equation}\label{equation-derivee}
(U-t)^2=3tU^2F(U)^2.
\end{equation}
%If we set $t=0$, we obtain $u^2$. 
By Theorem~\ref{thm-roots},
 this equation has two solutions, counted with multiplicities. Let us
 denote them  $U_1$ and $U_2$.
%Writing $U=\sum_i u_it^i$, and 
Using the first terms of $F(u)$,
one derives from~\eqref{equation-derivee} the first terms of the
series $U_i$. Remarkably,
one finds  $U_i= t +O(t^9)$ for $i=1,2$.

This observation leads us to \emm conjecture, that the series $U_i$ are the
same, so that~\eqref{equation-derivee} has a  only one solution,
of multiplicity
$2$. Let $\Delta(x_1,x_2,t,v)$ be 
the discriminant of $P(x_0,x_1,x_2,t,v)$ taken with respect to  $x_0$.
If $U_1=U_2\equiv U$, then, by Theorem~\ref{discriminant-original}, the
series $U$ is a root of 
$\Delta$ of multiplicity at least 4. As $\Delta$ factors as $tv^2D$,
for some polynomial $D\equiv D(v)$ of degree 
8  in $v$, our assumption
implies that for $0\leq i\leq 3$,
$$
\frac{\partial ^{i} D}{\partial v^i}(F(0), F'(0), t, U)=0.
$$
This gives  $4$ equations involving $3$ unknowns, namely $F(0)$,
$F'(0)$ and $U$.

Let us first eliminate  $F'(0)$ between $D(U)$ and
$D'(U)$. The resultant thus obtained reads
$t^{10}U^{6}(U-t)^6R_1$, where $R_1$ is a polynomial in $t$, $U$ and
$F(0)$, of degree
$8$ in $U$. The first few terms of $F(0)$ and $U$, which we have
computed,  rule out the possibility that $U=0$, but are not
sufficient to decide which  of the factors $(U-t)$ and $R_1$ are zero.

So let us first assume that $U=t$.
Taking the resultant in $F'(0)$ of $D(U)$
and $D''(U)$  gives  $F(0)=-t$. Returning to $\Delta$ provides
$F'(0)=1$. Set now $F(0)=-t$ and $F'(0)=1$ in
the original functional equation~\eqref{example-doubleroot}.
This gives the following cubic equation in $F(u)$:
$$
tu^2F(u)^3-(u-t)^2F(u)+(u-t)^3=0.
$$
By Theorem~\ref{thm-roots}, this equation has only one solution that is a \fps
\ in $t$. The form of the equation suggests to write $F=(u-t)G$, so
that $G$ satisfies $G=1+tu^2G^3.$

Hence, our assumption that~\eqref{equation-derivee} has a double root has led us to
the conjecture that the solution of~\eqref{example-doubleroot} is
$F=(u-t)G$, where $G$ 
is the unique series in $t$ satisfying $G=1+tu^2G^3$.
%$$F(t,u)=(u-t)(1+u^2t+3u^4t^2+12u^6t^3+\dots)$$
It is now straightforward to check that this series $F$ satisfies
$F(0)=-t$, $F'(0)=1$, and that the original functional
equation~\eqref{example-doubleroot} holds. 
Given that this equation has a unique power series solution, we have
solved it.   
% The fact that $F(u)$ is a multiple of $(u-t)$ confirms  that $U=t$
% is  a multiple root of~\eqref{equation-derivee}. 

%%%%%%%%%%%%%%%%%%%%%%%%%%%%%%%%%%%%%%%%%%%%%%%%%
\subsection{Applying the factorization theorem to the discriminant}
\label{section-factorization}
%%%%%%%%%%%%%%%%%%%%%%%%%%%%%%%%%%%%%%%%%%%%%%%%%
We exploit here the factorization theorem,
%   of Section~\ref{section-brown}, 
Theorem~\ref{thm-factorization}, in
combination with Theorem~\ref{discriminant-original}, which implies
%according to which 
that the discriminant 
$\Delta(F_1, \ldots, F_k,t,v)$ admits
% , in the generic case, 
$k$ multiple roots.

 Our example is again  the equation
for 3-constellations studied in Section~\ref{section-3const}. 
There, $k=2$, and the discriminant reads $\Delta(F_1,
F_2,t,v)=tv(v-1)^4R_1$, where $R_1$ is a polynomial in $F_1, F_2,
t$ and $v$, of degree $5$ in $v$. By
Theorem~\ref{discriminant-original}, $\Delta$ admits two
double roots $U_1$ and $U_2$, and we have seen 
% the first terms of the expansions of the $U_i$ show 
that they are actually double roots of $R_1$. What
about the fifth root of $R_1$? Setting $t=0$ in $R_1$ gives a
polynomial in $v$ of degree 4, so the fifth root of $R_1$ is infinite
at $t=0$. 

 Theorem~\ref{thm-factorization}, combined with the form of the series
 $U_i$, implies that $R_1$ factors as
\beq\label{factor-Delta-3const}
R_1=ct^p(1+tS)( 1+tvR) \Big( (1-v)^2+tvQ_0+tv^2Q_1\Big)^2
\eeq
where $S$, $R,Q_0$ and $Q_1$ belong to $\cs[[t]]$. Setting $t=0$
in this identity immediately gives $c=-4$ and 
$p=0$. Setting $v=0$ gives $S=0$.  Extracting the coefficient
of $v$ gives an expression of $Q_0$ in terms of $R$ and $F_1$. 
Extracting the coefficient
of $v^2$ gives an expression of $Q_1$ in terms of $R,F_1$ and $F_2$. 
% Similarly, the coefficient of $v^4$ provides an expression of
% $Q_1$. 
We now replace $Q_0$ and $Q_1$ by their expressions
in~\Ref{factor-Delta-3const}. 
The extraction of the coefficients of  $v^i$, for $i=3,4, 5$ gives a system of
3 polynomial equations relating $R, F_1$ and $F_2$. The elimination
of $R$ and $F_2$ yields back (after some heavy intermediate steps) the algebraic
equation~\eqref{3const-sol} satisfied by $F_1$. 

%Observe that we have not used the full force of the
%factorization~\eqref{factor-Delta-3const}. We could just as well have
%started from
%$$R_1=ct^p(1+tS)( 1+vW_3) \Big( 1+vS_0+v^2S_1\Big)^2.$$

 %%%%%%%%%%%%%%%%%%%%%%%%%%%%%%%%%%%%%%%%%%%%%%%%%
\subsection{Writing directly $k$ equations}
 %%%%%%%%%%%%%%%%%%%%%%%%%%%%%%%%%%%%%%%%%%%%%%%%%
We exploit here the results of Section~\ref{section-resultants}
(Theorem~\ref{2k-kcritere2}), which 
provide directly a system of $k$ equations between the series
$F_i$. Our guinea-pig is again  the equation 
for 3-constellations studied in Section~\ref{section-3const}. In order
to apply Theorem~\ref{2k-kcritere2}, we need two polynomials $P$ and
$Q$ having
$2$ roots in common. With the notation of
Section~\ref{section-3const}, these polynomials can be either $\Delta$
and $\Delta'_v$, or $R_1$ and $R_2$, or $R_3$ and $R_4$: the latter
pair being the simplest, we decide to start from it. Then $P$ has
degree $m=4$ and $Q$ has degree $n=8$. The Sylvester matrix of $P$ and
$Q$, denoted  $\mathcal
S_0(P,Q)$ in Section~\ref{section-resultants}, 
has size $12$.  By Theorem~\ref{2k-kcritere2}, its determinant $D_0$
--- the resultant 
of $P$ and $Q$ --- is zero, as well as the determinant $D_1$ of the matrix
$\mathcal S_1(P,Q)$, obtained by deleting the  last two columns,
as well as the last row of $a$'s (the 8th row) and the last row of
$b$'s (the 12th row). 

The  determinant $D_0$ is found to factor into two terms. The relevant
one has degree 8  in $F_1$ and  degree 4 in
$F_2$. The second determinant, $D_1$, does not factor, and has degrees
14 and 6 in $F_1$ and $F_2$ respectively. Still, Maple agrees to
eliminate $F_2$ between $D_1$ and the relevant factor of $D_0$. The
corresponding resultant contains four different factors, and the one that
vanishes  yields~\eqref{3const-sol} again. 

We observe that two of the above factors are squared. The occurrence
of repeated factors in iterated resultants is a systematic
phenomenon, which we discuss in Section~\ref{section-final}.

%%%%%%%%%%%%%%%%%%%%%%%%%%%%%%%%%%%%%%%%%%%%%%%%%
\section{The degree distribution of planar maps}
\label{section-distribution}
%%%%%%%%%%%%%%%%%%%%%%%%%%%%%%%%%%%%%%%%%%%%%%%%%
Let us return to the equations of Lemmas~\ref{lemma-finite-faces}
and~\ref{lemma-all-faces}, which 
characterize the face-distribution of rooted planar maps. In this
section, we solve these equations by
generalizing the approach of~\cite{bender-canfield}. 
Then we compare our solution to the result  obtained
in~\cite{BDG-planaires}  for the same problem.

%\subsection{A generalization of \cite{bender-canfield}}
%\label{general-bender-canfield}

\begin{Theorem}\label{generalization}
There exists a unique pair $(R_1,R_2)$ of \fps\ in $t$ with
coefficients in $\qs[z_1,z_2, \ldots]$ such that
\begin{equation}\label{R12}
R_1=\frac t2\sum_{i\ge1}z_i[u^{i-1}]R^{-1/2}
\quad \hbox{ and } \quad 
R_2 = t-3R_1^2 + \frac t2\sum_{i\ge1}z_i[u^i]R^{-1/2},
\end{equation}
where 
$$
R=1-4uR_1-4u^2R_2.
$$
 Let $G(t;\bm z)=G(t;z_1,z_2, \ldots)$ be the
generating function of  rooted planar maps, counted by the number
of edges (variable $t$) and the number of faces of degree $i$
(variable $z_i$). Then
$$
t^2 (tG(t;\bm z))'={(R_2+R_1^2)(R_2+9R_1^2)},
$$
where the derivative is taken with respect to $t$, and
$$
tG(t;\bm z)=
\frac 1 t \left( R_2+R_1^2
 \right)  \left( 3R_2+15R_1^{2}-2t \right) +
R_1  [u] \frac {\beta}{\sqrt R}- \frac 12 [u^2] \frac {\beta}{\sqrt R}
$$
where
$$\beta = \sum_{i\ge 1} z_i u^{-i}.$$
\end{Theorem}

\noindent{\bf Comments}\\
1. The equations defining $R_1$ and $R_2$ can also be written in terms of
$\beta$:
\beq\label{R-beta}
R_1=\frac t2 [u^{-1}]\frac {\beta}{\sqrt R} \quad \hbox{ and } \quad
R_2 = t-3R_1^2 + \frac t2[u^0]\frac {\beta}{\sqrt R}. 
\eeq
2. Let $m$ be a positive integer, and set
$z_i=0$ for $i >m$. Then $G(t;\bm z)$ is the face-distribution \gf\ of
planar maps in which all faces have degree at most $m$. The right-hand
sides of the equations defining $R_1$ and $R_2$ now involve only finitely
many terms, so that these two series are actually algebraic. The same
holds for $G(t;\bm z)$ (as stated in
Corollary~\ref{coro-face-algebraic}), and the above  
theorem makes this algebraicity explicit by providing a system of
three polynomial equations defining $R_1, R_2$ and $G$. For instance,
the \gf \ of planar maps in which all faces have degree 3, counted by
edges and faces, satisfies
$$
t^2G(t;z)={ \left( {R_2}+{{R_1}}^{2} \right)
  \left(3{R_2}+15R_1^{2}-2t- 24\,tz{R_1}{R_2} 
-56\,tz{{R_1}}^{3}
 \right) }
$$
with
$$
R_1={zt}  \left( { R_2}+3{{ R_1}}^{2} \right) \quad \hbox{ and
}\quad 
{ R_2}=t-3{{ R_1}}^{2}+  {zt}  \left( 6{ R_1}{ R_2}+10{{ R_1}}^{3}
 \right). $$

\smallskip
\noindent
{\bf Proof.} The existence and uniqueness of the series $R_j$ is
clear: think of extracting inductively from the equations~\eqref{R12} the
coefficient of $t^n$. The fact that $R_1$ and $R_2$ are multiples of $t$
implies that only finitely many values of $i$ are involved in this
 extraction, so that the coefficient of $t^n$ in $R_1$ and
$R_2$ is a \emm polynomial, in the $z_i$. 

We now want to relate  $R_1$
and $R_2$  to the face-distribution of planar maps. 
As noted in~\cite[p.~13]{bender-canfield}, it suffices to prove our results
when $z_i=0$ for $i >m$, for any $m\ge 3$.
%  (the series $F(u)$ then
%counts maps in which every finite face has degree at most $m$). 
Then  the equation of Lemma~\ref{lemma-finite-faces} may be written
the form $P(F(u), F_1, \ldots , F_{m-2},t,u)=0$:
\beq
\label{eqFcartes}
u^{m-2}F(u)=u^{m-2}+tu^mF(u)^2+t\theta_1(u)F(u)-t\sum_{j=0}^{m-2}u^{j}F_j\,
\theta_{j+2}(u),
\eeq
where $\theta_k(u) $ is the following polynomial in $u$,
of degree $m-k$:
$$
\theta_k(u)=\sum_{i=k}^mz_iu^{m-i}.
$$
Note that $F_0=1$. This equation coincides with  Eq.~(2.2)
of~\cite{bender-canfield}, apart for the  value of $\theta_k$.

\noindent We now apply the general strategy of Section~\ref{section-key}. 
The condition $P'_0(F(U), F_1, \ldots , t,U)=0$ reads:
$$
U^{m-2}=2tU^mF(U)+t\theta_1(U).
$$
By Theorem~\ref{thm-roots}, this equation has $m-2$ solutions, $U_1, \ldots ,
U_{m-2}$, which are fractional series in $t$ (with coefficients in an
algebraic closure of $\qs(z_1, \ldots , z_m)$). All of them vanish at
$t=0$. By Theorem~\ref{discriminant-original}, these series 
are multiple roots of the discriminant
$\Delta(u)\equiv \Delta(F_1, \ldots , F_{m-2},t,u)$. This discriminant
is found to be
\beq
\label{Delta-face}
\Delta(u)=\left(t\theta_1(u)-u^{m-2}\right)^2-4tu^m\left(u^{m-2}-t\sum_{j=0}^{m-2}u^{j}F_j\theta_{j+2}(u)\right).
\eeq
It has degree (at most) $2m-2$ in $u$, and it reduces to $u^{2(m-2)}$
when $t=0$. By 
the Newton-Puiseux theorem, this implies that the series $U_i$, for
$1\le i \le m-2$,  are the
only roots of $\Delta(u)$ that are finite at $t=0$, and that they have
multiplicity 2 exactly. The remaining 
roots are infinite. The canonical factorization of $\Delta(u)$
(Theorem~\ref{thm-factorization}) thus reads:
$$
\Delta(u)= c t^p (1+tS(t)) \Big( 1+tuS_1(t,u)\Big) \Big(u^{m-2}+
tS_2(t,u)\Big)^2
$$
where $S_1$ has degree (at most) 1 in $u$ and $S_2$ has degree  at
most $m-3$.
Setting $t=0$ in~\eqref{Delta-face} shows that $p=0$ and $c=1$. Setting
$u=0$ then gives 
$t^2(1+tS(t))S_2(t,0)^2=t^2z_m^2$, and we finally choose to
write the canonical factorization of $\Delta(u)$ with the notation
of~\cite{bender-canfield}:    
$$
\Delta=RQ^2
$$
with
$$
R=1-4uR_1-4u^2R_2  \quad \hbox{ and }   \quad Q=tz_m+\sum_{i=1}^{m-2}Q_iu^i.
$$
 The $R_i$ and $Q_i$ are power series in
$t$ with coefficients in $\qs(z_1, \ldots , z_m)$.

The derivations of the equations defining the $R_i$, and of the
expression of $(tG)'$, now faithfully
follow~\cite{bender-canfield}. Let us simply recall where \eqref{R12}
 comes from. Solving~\eqref{eqFcartes} gives
\begin{equation}\label{F-racine-quadratique}
2tu^mF(u)=u^{m-2}-t\theta_1(u)\pm\sqrt{\Delta(u)}=
u^{m-2}-t\theta_1(u)+Q\sqrt R,
\end{equation}
so that
\begin{equation}\label{congruence-Q}
Q=\frac{2tu^m+t\theta_1(u)-u^{m-2}}{\sqrt R} +O(u^{m+1}).
\end{equation}
Recall that $Q(u)$ has degree $m-2$  in $u$. Extracting the
coefficients of $u^{m-1}$ and $u^m$  in the above identity gives
\eqref{R12}.

Let us finally derive an expression for $G(t;\bm z)$. By~\eqref{congruence-Q},
\beq
\label{Qi}
Q_i= [u^i]\frac{t\theta_1(u)-u^{m-2}}{\sqrt R} \quad \hbox{ for } 0\le
i \le m-2.
\eeq
Now by Lemma~\ref{lemma-all-faces} and~\eqref{F-racine-quadratique},
$$
\begin{array}{lllll}
2t^2G= 2t{F_2}  &= &[u^{m+2}] Q\sqrt R= \displaystyle \sum_{i=0}^{m-2} Q_i[u^{m+2-i}]
\sqrt R\\
&=&\displaystyle \sum_{i=0}^{m-2} [u^i]\frac{t\theta_1(u)-u^{m-2}}{\sqrt R}[u^{m+2-i}]
\sqrt R \hskip 30mm
\hbox{by \eqref{Qi}} \\
&=&\displaystyle \sum_{i=0}^{m+2} [u^i]\frac{t\theta_1(u)-u^{m-2}}{\sqrt R}[u^{m+2-i}]
\sqrt R- \sum_{i=m-1}^{m+2} [u^i]\frac{t\theta_1(u)-u^{m-2}}{\sqrt
  R}[u^{m+2-i}] \sqrt R\\
& =&\displaystyle[u^{m+2}]\left(t\theta_1(u)-u^{m-2} \right)- \sum_{i=m-1}^{m+2} [u^i]\frac{t\theta_1(u)-u^{m-2}}{\sqrt
  R}[u^{m+2-i}] \sqrt R\\
%&=&-\displaystyle  \sum_{i=m-1}^{m+2} [u^i]\frac{t\theta_1(u)-u^{m-2}}{\sqrt
%  R}[u^{m+2-i}] \sqrt R\\
&=&-\displaystyle  \sum_{i=-1}^{2} [u^i]\frac{t\beta-u^{-2}}{\sqrt
  R}[u^{2-i}] \sqrt R,
\end{array}
$$
where $\beta=u^{-m}\theta_1(u)$. The expected expression of $G(t;\bm
z)$ follows, upon using~\eqref{R-beta}. 
\cqfd

%\subsection{Comparizon with \cite{BDG-planaires}}
In \cite{BDG-planaires}, another characterization of the
face-distribution of planar maps was obtained, using two different
methods: first, using matrix integrals, and then by a purely
combinatorial approach. Both methods yield the same expression for the
series $G(t;\bm z)$, but this expression differs from that of
Theorem~\ref{generalization}. Our aim is now to relate 
these two different expressions. Let us first recall the expression
of~\cite{BDG-planaires}. 

%The expression of~\cite{BDG-planaires} also involves a pair of
%auxilliary series, $S_1$ and $S_2$. 
%
\begin{Theorem}[\cite{BDG-planaires}]\label{bdg-thm}
There exists a unique pair $(S_1,S_2)$ of \fps\ in $t$ with
coefficients in $\qs[z_1,z_2, \ldots]$ such that
\begin{equation}\label{S12}
S_1 = t[v^0]W
\quad \hbox{ and } \quad S_2 = t+t[v^{-1}]W
\end{equation}
where 
$$
W=\sum_{i\geq1}{z_i}P^{i-1} \quad \hbox{ and } \quad P=v+S_1+S_2/v.
$$
The face-distribution 
generating function of  rooted planar maps, denoted above $G(t;\bm
z)$, satisfies
$$
tG(t;\bm z)=S_1^2+S_2-2S_1[v^{-2}]W-[v^{-3}]W.
$$
\end{Theorem}

\begin{Propo}\label{comparizon}
The solutions to the face-distribution  problem given by
Theorems~{\rm \ref{generalization}} and {\rm \ref{bdg-thm}} are
related as follows.
\begin{itemize}
\item[$(i)$]  The auxilliary series $R_i$ and $S_i$ satisfy
$$
S_1=2R_1 \quad \hbox{ and } \quad  S_2=R_2+R_1^2.
$$
Moreover, for all $\ell\ge 0$,
$$
[u^\ell] \frac \beta{\sqrt R} =[v^0] P^{\ell+1} W.
$$
\item[$(ii)$] The following identities, valid for all $k\ge 0$ and $j\in \zs$,
\begin{eqnarray*}
[v^j]P^{k+1}W&=& [v^{j-1}]P^{k}W+ S_1 [v^{j}]P^{k}W+ S_2
[v^{j+1}]P^{k}W, \\
\left[\right.\! v^j\!\left.\right]P^kW&=&  S_2^{-j}[v^{-j}]P^kW,
\end{eqnarray*}
allow one to express any term $[v^j]P^kW$ as a linear combination
of terms $[v^{-i}]W$, for $i\ge 0$, with coefficients in
$\qs(S_1,S_2)$.\\
\item[$(iii)$] Rewrite the expression of $G(t;\bm z)$ given in
Theorem~{\rm\ref{generalization}} in terms of $S_1$, $S_2$ and $P$ using
$(i)$. Then, use $(ii)$ to rewrite this in terms of  $S_1$, $S_2$ and the
$[v^{-i}]W$, for $i\ge 0$. Use finally~{\rm\eqref{S12}} to express
$[v^{0}]W$ and $[v^{-1}]W$ in terms of $S_1$ and $S_2$: 
The resulting expression of $G(t;\bm z)$ is that of
Theorem~{\rm\ref{bdg-thm}}.
\end{itemize}
\end{Propo}

\noindent
{\bf Proof. }
$(i)$ Let us introduce the series 
$$
\bar R_1:= \frac{S_1}2, \quad  \bar R_2:= S_2-\frac{S_1^2} 4 \quad
\hbox{and} \quad \bar R:= 1-4u\bar R_1 -4u^2 \bar R_2 =
(1-uS_1)^2-4u^2S_2.
$$
We want to prove that $\bar R, \bar R_1$ and $\bar R_2$
satisfy~(\ref{R12}) (with bars  over all unknowns).
In view of~\eqref{S12}, the first equation in~\eqref{R12} holds if and only if
$$
[v^0]W = \sum_{j\geq 0}z_{j+1}[u^{j}]\bar R^{-1/2}.
$$
Given that $W=\sum_j z_{j+1} P^{j}$, it suffices to prove that for all
$j\ge 0$, 
\beq\label{id-base}
[v^0]P^j= [u^j] \bar R^{-1/2},
\eeq
or, upon taking \gfs, that
$$
\sum_{j\ge 0} u^j [v^0]P^j=  \bar R^{-1/2}.
$$
But 
\begin{eqnarray*}
\sum_{j\ge 0} u^j [v^0]P^j&=&
\sum_{j\ge 0}u^j\sum_{k\ge 0} {j  \choose {2k}}  {{2k} \choose k}
S_1^{j-2k}S_2^k\\
% & =& \sum_{k\ge 0} {{2k} \choose k}
%S_2^k\sum_{j\ge 0}{j  \choose {2k}}u^jS_1^{j-2k}\\
&=&\sum_{k\ge 0} {{2k} \choose k} S_2^k
\frac{u^{2k}}{(1-uS_1)^{2k+1}}\\ 
&=& \frac 1 {1-uS_1} \left(1-\, \frac{4u^2S_2}{(1-u
  S_1)^{2}}\right) ^{-1/2} = \bar R^{-1/2}.
\end{eqnarray*}
(By convention, ${b \choose a} =0$ unless $0\le a\le b$.)
The first equation of~\eqref{R12} follows. The second one reads, in
view of~\eqref{S12},
$$
2[v^{-1}]W+S_1[v^0]W= \sum_{i\geq 1}z_i[u^{i}]\bar R^{-1/2}.
$$
In order to prove it, it suffices to check that for all $j\ge 0$,
\begin{eqnarray*}
2[v^{-1}]P^j+S_1[v^0]P^j&=&[u^{j+1}]\bar R^{-1/2},\\
&=&[v^0]P^{j+1} \hskip 20mm \hbox{ by \eqref{id-base}}.
\end{eqnarray*}
This is easily proved by first extracting the coefficient of $v^0$ in
$P^{j+1}= (v+S_1+S_2/v)P^j$, and then noticing that
$S_2[v]P^j=[v^{-1}]P^j$ (this comes from the fact that $P^j$ is left
unchanged when replacing $v$ by $S_2/v$). 
Since $\bar R_1$ and $\bar R_2$ satisfy~\eqref{R12}, they 
coincide respectively with the series $R_1$ and $R_2$.
The second result of $(i)$ follows from~\eqref{id-base}.

 The first identity of $(ii)$ is simply obtained by writing
$$
P^{k+1} W= (v+S_1+S_2/v)P^kW,$$
and extracting the coefficient of $v^j$. The second one follows from
the fact that $P$, and hence $W$, is left invariant upon replacing $v$
by $S_2/v$.

Finally, $(iii)$ is a straightforward verification.
 \cqfd

%%%%%%%%%%%%%%%%%%%%%%%%%%%%%%%%%%%%%%%%%%%%%%%%%
\section{Hard particles on planar maps}
\label{section-hard}
%%%%%%%%%%%%%%%%%%%%%%%%%%%%%%%%%%%%%%%%%%%%%%%%%
Let us return to the equations established in
Lemma~\ref{lemma-hard-particles} for planar maps carrying
hard particles. We will solve this system when $x=y=1$. That is, the series
$F(u)\equiv F(t,s,u)$ counts maps rooted at a vacant vertex by the total number of
edges (variable $t$), the number of frustrated edges (variable $s$)
and the number of white corners in the root-face (variable $u$). The
series $G(u)\equiv G(t,s,u)$ counts maps rooted at an occupied vertex, according to
the same statistics. 

The first step consists in 
%reducing the system to a single equation, obtained by
 eliminating $G(u)$.
% between the two equations. 
This
gives an  equation of the form $P(F(u),F(1),G(1),t,u)=0$, which is
cubic in $F(u)$.   

The next steps require a computer, but otherwise copy faithfully the
bare-hands strategy of Section~\ref{section-3const}. We do not give
the details. Let us simply mention that, when $s=1$, the two series
$U_i$ that cancel $P'_0(F(u),F(1),G(1),t,u)$ are formal power series
in $\sqrt t$:
$$
U_{1,2}=1+{t}\pm{t}^{3/2}+4{t}^{2}\pm 17/2{t}^{5/2}+O(t^3)
$$
while, when $s\not =1$, they are formal series in $t$ with
coefficients in $\qs(s)$:
$$
U_1=1+t+{\frac { \left( 3\,s-4 \right) {t}^{2}}{s-1}}+{\frac { \left( -25+
64\,s+{s}^{4}+12\,{s}^{3}-51\,{s}^{2} \right) {t}^{3}}{ \left( s-1
 \right) ^{3}}} +O(t^4),
$$
$$
U_2=1+st+{\frac {{s}^{2} \left( -2+3\,s \right) {t}^{2}}{s-1}}+{\frac {{s}
^{2} \left( -28\,{s}^{3}+13\,{s}^{4}+16\,{s}^{2}-2 \right) {t}^{3}}{
 \left( s-1 \right) ^{3}}} +O(t^4).
$$
Another interesting observation is that, in all cases, some of the
roots of the  discriminant $\Delta(u)\equiv\Delta(F(1),G(1),t,u)$ that are
finite at $t=0$ are \emm simple,.  For
instance, when $s=1$, this discriminant  has two simple roots, $U_3$
and $U_4$, of the following form:
$$
U_{3,4}=1\pm2 \,i{t}^{1/2} -5\,{t}\mp 13\,i{t}^{3/2}+O(t^2).
$$
In other words, $\Delta(u)$ does not have  a square root in
$\GK[u][[t]]$ (by Proposition~\ref{propo-even-mult}).  The other roots of
$\Delta(u)$ that are finite at $t=0$ are $U_1$, $U_2$, $0$ and $1$,
and have an even multiplicity.    
% This confirms the result already announced at the end of
%Section~\ref{section-brown}: it is not always true that the the roots of the
%discriminant are multiple as soon as they are finite at $t=0$. 

At the end of the elimination procedure, one obtains a pair of quartic
polynomial equations for the 
series $F(1)$ and $G(1)$. The corresponding curves have genus $0$ (as
is often the case  for hard-particle models~\cite{BDFG02a}), and we have
found, with the help of the {\tt algcurves} package of {\sc Maple}, a simple
parametrization of them. In this form, our results are begging for a
purely combinatorial derivation, in the vein
of~\cite{mbm-schaeffer-ising,bdg-pdures,bdg-pdures-bip}. 
\begin{Propo}
Let $T\equiv T(t,s)$ be the unique \fps\ in $t$ with constant term $0$
satisfying
$$
T(1-2T)(s-3T+3T^2)=s^2t.
$$
Then $T$ has actually coefficients in $\ns[s]$. Morever, the \gfs\
$F(t,s,1)\equiv F(1)$ and $G(t,s,1)\equiv G(1)$ that count  planar
maps carrying hard particles (rooted, respectively,  at 
an empty and an occupied vertex) satisfy
\begin{eqnarray*}
s^3t^2F(1)&=& {T}^{2}
 \left(s -4\,T-3\,sT+15\,{T}^{2}+s{T}^{2}-15\,{T}^{3}+4\,{T}^
{4} \right), \\
s^4t^2G(1)&=&  {T}^{3}\left( s-3\,T+3\,{T}^{2} \right) 
 \left(s -4\,T-3\,sT+14\,{T}^{2}-9\,{T}^{3} \right) .
 \end{eqnarray*}
\end{Propo}
\noindent
{\bf Proof.}
The only result that does not follow from the
elimination of $U_1$ and $U_2$ is the fact that $T$ has coefficients
in $\ns[s]$. By writing 
$$
T=sS,\quad S= \frac{t}{(1-2sS)(1-3S(1-sS))},\quad
S(1-sS) = t\, \frac{1-sS}{1-2sS}\, \frac 1{1-3S(1-sS)},
$$
it is easy to prove by induction on $n $ that the coefficient of $t^n$
in $S$ and in $S(1-sS)$ belongs to $\ns[s]$. 
Since $T=sS$, the same holds for the
coefficients of $T$.
\cqfd

By studying the singular behaviour of the series $F(1)$ and $G(1)$
when $s=1$, one obtains the following corollary.
\begin{Coro}
The number of $n$-edge planar maps carrying hard particles 
is equivalent to
$$
\alpha \left( \frac{\sqrt [3]{39509+23436\,\sqrt {62}}-\sqrt [3]{-39509+23436\,\sqrt {62}}
+38}3
\right) ^n n^{-5/2} \simeq \alpha (15.4...)^n n^{-5/2},
$$
for some positive constant $\alpha$.
\end{Coro}
\noindent {\bf Proof.} We apply the general principles that relate the
singularities of an algebraic series to the asymptotic behaviour of
its coefficients~\cite{fla-poly-alg}. 
%  The argument is the same for the series $F(1)$ and $G(1)$. 
Since  $T$, $F(1)$ and $G(1)$ have
non-negative coefficients, their radius of convergence is one of their
singularities. The expressions of $F(1) $ and
$G(1)$ show that their singularities are also
singularities of $T$. As the leading coefficient of the equation
defining $T$ does not vanish, the singularities of $T$ are among 
the roots of the discriminant of its minimal polynomial, that is,
among the roots of
$$
\delta=18432t^3-1545t^2+38t-1.
$$
The only real root of $\delta$ is $\rho\simeq.065$.  Hence
$\rho$ is the radius of convergence of $T$, $F(1)$ and $G(1)$. The
modulus of the other two roots of $\delta$  is less than $\rho$.  So
$\rho$ is actually the only  singularity 
of $T$, $F(1)$ and $G(1)$. A local expansion of these three series  in
the neighborhood of $t=\rho$ shows that $T$ has a square root type
singularity, while $F(1)$ and $G(1)$ have a singularity in
$(1-t/\rho)^{3/2}$.
 This implies
that the $n$th coefficient of $F(1)$ and $G(1)$ is asymptotic to $\alpha
\rho^{-n} n^{-5/2}$ for some positive constant $\alpha$ (which is not
the same for $F(1)$ and $G(1)$). But
$\rho^{-1}$ is exactly the constant occurring above in the  corollary.
\cqfd

\noindent{\bf Note.} A similar study can be conducted for a generic
value of $s\in (0, +\infty)$. We have not worked out all the details,
but it seems that the above pattern persists for all $s\in \rs_+$. In
other words, there is no (physical) phase transition in this model. At
any rate, it is not very hard to prove that the radius of convergence
$\rho(s)$ is a smooth function of $s$, equal to the branch of
$$
18432\,{s}^{4}{\rho}^{3}-3\,{s}^{2}  {\rho}^{2}\left( 963-2496\,s+2048\,{s}^{2}
 \right)+2\,\rho \left( 16\,{s}^{2}+21\,s-18 \right)  \left( 4
\,s-3 \right) ^{2}- \left( 4\,s -3\right) ^{3}
$$
that equals $1/12$ at $s=0$.
%%%%%%%%%%%%%%%%%%%%%%%%%%%%%%%%%%%%%%%%%%%%%%%%%%
\section{Concluding remarks and questions}
\label{section-final}
Let us begin by bragging a bit about some positive points of this
paper. We have proved that the series $F(t,u)$ given by a functional
equation of a certain type (see~\eqref{main-eq}) are  algebraic. As
illustrated in Section~\ref{section-equations}, this tells us that a number of
enumerative problems have an algebraic \gf , without having to solve
them in detail. Our general strategy gives a system of $3k$ polynomial
equations. Its reduction to $2k$ equations
(Section~\ref{section-discriminant}) has a theoretical 
interest, and tells us what is left of the ``quadratic method'' for
equations that are no longer quadratic.

However, the practical aspects of our approach probably require more
work. Generally speaking, we are lacking an efficient elimination
theory for polynomial systems which, as~(\ref{key-P}--\ref{key-Pu})
or~\eqref{system-2k}, are 
highly symmetric. The case of $3$-constellations, solved by different
approaches in Section~\ref{section-practise}, shows that even when the final result is
relatively simple (here, $F_1$ satisfies a cubic equation), the
intermediate steps may involve  big polynomials. Does this mean
that 3-constellations are somehow pathological, or that we have not
conducted the calculations in the best possible way (or both...)? We
have no definite answer to this question, but the following
observations may be of some interest:

\noindent {\bf 1.}
The degree of $F(u)$ may be very big compared to the degrees
  of the
  original functional equation. Consider, for instance, the
  enumeration of walks on the half-line $\ns$, that start 
  from $0$ and take steps $+k$ and $-\ell$, where $k$ and $\ell$ are
  coprime. The case $(k,\ell)=(3,2)$ 
  was solved in Section~\ref{section-kernel}. In general, the equation reads
$$
(u^\ell-t(1+u^{k+\ell}))F(u)=u^\ell -t \sum_{i=0}^{\ell-1} u^i F_i.
$$
Its solution satisfies~\cite[Ex.~4]{bousquet-petkovsek-1}:
$$
tF_0=tF(0)=(-1)^{\ell+1} \prod_{i=1}^\ell U_i,
$$
where $U_1, \ldots , U_\ell$ are the $\ell$ roots of the kernel that
are finite at $t=0$. It 
can be proved that $F_0$ has degree exactly ${{k+\ell }\choose
  k}$. When $\ell=k-1$, this degree is exponential in $k$, even though
the original equation is linear in $F(u)$ and all the $F_i$ (and has
degree $2k-1$ in $u$).

\smallskip
\noindent {\bf 2.} Certain resultant calculations yield systematically repeated
   factors. Imagine we are trying to find a polynomial equation for
   $F_1$, starting from $P(F(u),F_1,F_2,t,u)=0$. At some point, we end
   up with two polynomials $R(u)$ and $Q(u)$, with coefficients in
   $\GK[F_1,F_2,t]$, which have two roots in common. We thus apply
   Theorem~\ref{2k-kcritere2}: the determinants $D_0$ and $D_1$ of
   $\mathcal S_0$ and    $\mathcal S_1$  vanish. We take the
   resultant of $D_0$ and $D_1$ in $F_2$ to obtain a polynomial
   equation for $F_1$. 
%It is known that
% , if $R$ and $Q$ were \emm    generic, polynomials in $u$, $F_1$ and
%   $F_2$ of a total degree $m$    and $n$,  
Then every factor in this resultant has multiplicity
   at least 2~\cite[Thm.~3.4]{mccallum}. 

Moreover, a similar reduction might apply to equations with a single
unknown function $F_1$. For such equations, we obtain a polynomial
equation for $F_1$ by writing that  the iterated discriminant
$discrim_v (discrim_x P(x,F_1,t,v))$ vanishes. Again, if $P$ is a
generic polynomial in $x, F_1$ and $v$ of total degree $d\ge 3$, then it is
conjectured that this iterated discriminant has repeated
factors~\cite[p.~384]{mccallum}. Note that this does not mean that
we will always meet such a factor in our examples, since they do not have
generic coefficients. This is illustrated by the following
example, 
which is a generalization of the equation for planar maps
(Section~\ref{section-quadratic}):
$$
F(u)  =1+atu^2F(u)^2+ tu \, \frac{uF(u)-F(1)}{u-1}.
$$
There, the iterated discriminant is an irreducible polynomial of
degree 4 in $F_1$ -- but the equation we start from is  not
a   \emm generic, equation of total degree 5.

\medskip
Finally, let us underline that it is still an unsolved problem to
enumerate $m$-constellations starting from the equations of
Proposition~\ref{propo-constellations} (even though the result is known
to be remarkably simple~\cite{mbm-schaeffer-constellations}). 

%%%%%%%%%%%%%%%%%%%%%%%%%%%%%%%%%%%%%%%%%%%%%%
\bigskip
\noindent{\bf Acknowledgements.} We thank Pierrette
Cassou-Noguès,  Emmanuel Briand,
 Frédéric Chyzak, Bruno Salvy and Alin
Bostan for interesting discussions and references.
We are also grateful to Gilles Schaeffer for communicating the
equations of Lemma~\ref{lemma-hard-particles}.

\bibliographystyle{plain}
\bibliography{biblio.bib}

\end{document}